\documentclass{amsart}

\newtheorem{thm}{Theorem}[section]
\newtheorem{prop}[thm]{Proposition}
\newtheorem{corollary}[thm]{Corollary}
\newtheorem{lemma}[thm]{Lemma}

\theoremstyle{definition}
\newtheorem{dfn}[thm]{Definition}
\newtheorem{ex}[thm]{Example}

\theoremstyle{remark}
\newtheorem{rem}[thm]{Remark}

\numberwithin{equation}{section}

\newcommand{\grad}{\operatorname{grad}}

\newcommand{\id}{\operatorname{id}}

\newcommand{\F}{\ensuremath{\mathcal{F}}}
\newcommand{\CP}{\ensuremath{\mathcal{P}}}
\newcommand{\CJ}{\ensuremath{\mathcal{J}}}
\newcommand{\CC}{\ensuremath{\mathcal{C}}}
\newcommand{\CD}{\ensuremath{\mathcal{D}}}
\newcommand{\CG}{\ensuremath{\mathcal{G}}}
\newcommand{\CH}{\ensuremath{\mathcal{H}}}

\newcommand{\RR}{\mathbb R}
\newcommand{\ZZ}{\mathbb Z}

\newcommand{\prf}{{\it Proof.\, }}
\newcommand{\eop}{\hfill $\Box$\medskip}
\newcommand{\eps}{\varepsilon}

\newcommand{\Verweist}{Theorem\ }
\newcommand{\Verweisp}{Proposition\ }
\newcommand{\Verweisl}{Lemma\ }
\newcommand{\Verweisc}{Corollary\ }
\newcommand{\Verweisr}{Remark\ }

\begin{document}
\title{Parallel Focal Structure and Singular Riemannian Foliations}
\begin{abstract}
We give a necessary and sufficient condition for a submanifold with parallel focal structure to give rise to a 
global foliation of the ambient space by parallel and focal manifolds. We show that this is a singular Riemannian foliation with complete orthogonal transversals. For this object we construct an action on the transversals that generalizes the Weyl group action for polar actions.
\end{abstract}
\author{Dirk T\"oben}
\address{Mathematisches Institut, Universit\"at zu K\"oln, Weyerthal 86-90, 50931 K\"oln, Germany}
\email{dtoeben@math.uni-koeln.de}
\date{June 11, 2004}
\maketitle
\section{Introduction}
In his thesis \cite{ewert}, Ewert introduced the notion of a submanifold with parallel focal structure as a generalization of isoparametric submanifolds in Euclidean space (see for instance \cite{palaisterng}) and equifocal submanifolds in simply connected symmetric spaces (\cite{terngthorbergsson}). For a survey on these and related objects, see \cite{thorbergssonsurvey}. In order to generalize results known for isoparametric and equifocal submanifolds we pursue a completely different approach than in the respective theories. We first look for a condition under which the given submanifold with parallel focal structure gives rise to a singular foliation of the ambient space. Then we will obtain similar properties as for isoparametric and equifocal submanifolds as a consequence of this foliated structure.\par  
Let $N$ be a complete Riemannian manifold and $M$ be a submanifold. We begin by asking under which conditions we have a good partition of $N$ by parallel and focal manifolds of $M$. To properly define parallel and focal manifolds, $M$ has to suffice the following minimal conditions.
\begin{enumerate}
\item $\nu M$ is flat, i.e., any $v\in\nu M$ can be locally extended to a parallel normal field $v'$.
\item the rank of the locally defined map $\exp\circ\, v'$ is constant for any $v\in\nu M$.
\end{enumerate}
For $v\in\nu M$ we define 
$$
M_v=\left\{\exp \Big(({\overset{1}{\underset{0}\parallel}}c)v\Big)\ \bigg|\ c\  \mbox{is a curve in}\ M \right\},
$$ 
where $\parallel c$ denotes normal parallel translation along $c$. We call $M_v$ {\it parallel manifold} of $M$ if the rank of $\exp\circ\, v'$ is maximal, otherwise {\it focal manifold}.
\begin{dfn}(\cite{heintze})
We say that $M$ gives rise to a {\it global 
foliation} $\F=\{M_v\ |\ v\in\nu_pM\}$ of $N$,  if $\bigcup \F=N$, and $M_v\cap M_w\neq\emptyset$ implies $M_v=M_w$.
\end{dfn}
\begin{ex}\label{bspfol}
First we take $N=S^2$ and $M$ a parallel of the equator. Clearly $M$ induces a global foliation. Next we consider the flat torus $N=T^2$ and a small distance circle $M$ centered at a point $p$ in $N$. $M$ does not induce a global foliation of $N$. 
\end{ex}
With the notion of the cut locus of a proper immersion we give a necessary condition for $M$ under the above minimal conditions to induce a global foliation in \Verweisp\ref{necessary_foliation}. This condition is not sufficient. Now assume that $M$ fulfills the minimal conditions and in addition that through any $x\in M$ there is a complete totally geodesic submanifold $\Sigma_x$, a {\it section}, with $T_x\Sigma_x=\nu_xM$. This class of submanifolds generalizes isoparametric submanifolds in Euclidean space and equifocal submanifolds in simply connected symmetric spaces for arbitrary ambient spaces. For this class we prove that the condition, that any two sections through a regular value of the normal exponential map of $M$ coincide, is necessary and sufficient to give rise to a global foliation. We then call $M$ a {\it submanifold with parallel focal structure} (the precise definition is given in \ref{defsub}). We only formulate the sufficiency condition:\smallskip

{\noindent \bf Theorem.} {\it A closed and embedded submanifold $M$ with parallel focal structure and finite normal holonomy of a complete Riemannian manifold gives rise to a global foliation $\F$. This global foliation is a singular Riemannian foliation admitting sections.}\smallskip

The definition of a singular Riemannian foliation is given below. The converse statement, that a regular leaf of a singular Riemannian foliation admitting sections has parallel focal structure was proven by Alexandrino in \cite{alexandrino}. 
\begin{dfn}[\cite{molino}]
Let $\F$ be a partition of injectively immersed submanifolds (the {\it leaves}) of a Riemannian manifold $N$. For any $p\in N$ let $M_p$ be the leaf through $p$ and let $T\F=\bigcup_{p\in N}T_pM_p$. We define $\Xi(\F)$ as the module of (differentiable) vector fields on $N$ with values in $T\F$. We call $\F$ a {\it singular Riemannian foliation}, if
\begin{enumerate}
\item (Transnormality) a geodesic starting orthogonally to a leaf intersects the leaves it meets orthogonally;
\item (Differentiability) $\Xi(\F)$ acts transitively on $T\F$, i.e., for any $v\in T_p\F, p\in N$ there is $X\in\Xi(\F)$ with $X_p=v$.
\end{enumerate}
A leaf of maximal dimension is called {\it regular}, and so each point of it, otherwise {\it singular}. If, in addition, for any regular $p$ there is an isometrically immersed complete totally geodesic submanifold $\Sigma_p$ (the {\it section}) with $T_p\Sigma=\nu_pM_p$, that meets any leaf and always orthogonally, $\F$ is a singular Riemannian foliation {\it admitting sections}.
\end{dfn}
A partition into injectively immersed submanifolds with every leaf of the same dimension is a foliation if and only if (2) is fulfilled. A foliation with (1) is a Riemannian foliation. The set of orbits of an isometric Lie group action on a Riemannian manifold $N$ is a singular Riemannian foliation. The set of orbits of a polar action is a singular Riemannian foliation admitting sections.\medskip

In section 2 we introduce the cut locus of a submanifold.\par
In section 3 we prove the above theorem. In 3.1 we describe general properties of submanifolds with parallel focal structure. In 3.2 we associate a regular Riemannian foliation $(\hat N,\hat\F)$ to $M$, the {\it blow-up}, which was constructed in \cite{boualem} for singular Riemannian foliations admitting sections (with relatively compact leaves). With this blow-up we show that each parallel manifold of $M$ also has parallel focal structure and we derive the theorem in 3.3.\par
In section 4 we study singular Riemannian foliations. In 4.1 we give an alternative proof of the converse of the above theorem. In 4.2 we introduce an action on the sections, the {\it transversal holonomy group}, which generalizes the Weyl group action of a polar action and we give applications. Furthermore we prove that each regular leaf of a singular Riemannian foliation admitting sections in a simply connected symmetric space has trivial normal holonomy.\par
The author would like to thank G. Thorbergsson for many helpful discussions.

\section{Cut Locus of a Proper Immersion}
In this section we introduce the notion of a cut locus of a submanifold. This is a generalization of the cut locus of a point which is defined for instance in \cite{klingenberg}. With the aid of this notion we can formulate a necessary condition for a properly immersed submanifold with the two minimal conditions stated in the first section to give rise to a global foliation.\par
Let $N$ be a complete and connected Riemannian manifold and $M$ a manifold. By $\varphi:M\to N$ we will always denote an isometric immersion. Let $\iota:\nu M\to TN$ be the canonical immersion. We write $\eta:=\exp^\perp=\exp\circ\,\iota:\nu M\to N$ and $\eta^r:B_r(\nu M) \to N$ for the restriction of $\eta$ to
the normal ball bundle of $M$ of radius $r$. Let $\varphi:M\to N$ be a proper immersion. Let $\gamma_v$ denote the geodesic with initial
vector $v\in TN$.\par 
\begin{dfn} 
We define $\sigma:\nu^1 M\to [0,\infty]$ by
$$
\sigma(v)=\sup\{t\in \RR\ |\ d(\gamma_{\iota(v)}(t),\varphi(M))=t\}.
$$
We call $\sigma(v)$ the {\it cut distance of $M$ in direction
$v$}. The {\it normal cut locus} ${\CC}^{\nu M}$ of $\varphi$ is defined by ${\CC}^{\nu M}:= \{\sigma(v)v\ |\
\sigma(v)<\infty, v\in\nu^1M\}$ and the {\it cut locus} ${\CC}_\varphi$ or $\CC_{(M,N)}$ by $\exp^\perp {\CC}^{\nu M}$.
\end{dfn}

\begin{dfn} 
A vector $v\in\nu M$ is called a {\it normal} or {\it normal vector} and the geodesic $\gamma_{\iota(v)}$ a {\it normal geodesic}. If $\|v\|\leq \sigma(v/\|v\|)$ for  $v\in\nu M$, we call $v$ {\it minimal} and $\gamma_v|[0,1]$, and any reparametrization of constant speed, a {\it minimal geodesic (segment)}. This terminology is justified by the fact that, in the set $\eta^{-1}(p)$, the minimal vectors have the least length.
We call a normal vector $v$ a {\it focal normal} and $\eta(v)$ a {\it focal point}, if $v$ is singular with respect to $\eta$.
We call a minimal vector $v$ a {\it cut vector} and $\eta(v)$ a {\it cut point}, if there is a minimal $w\in\nu M$ with $\iota(w)\neq\iota(v)$ having the same endpoint as $v$. In this case $\|v\|=\sigma(v/\|v\|)=\sigma(w/\|w\|)=\|w\|$. 
\end{dfn}

It is easy to see that the limit of a converging sequence of minimal normal vectors is minimal and it is kown that $tv$ for $t>1$ is not minimal if $v$ is a focal vector or a cut vector. Also note that for every $p\in N$ there is a shortest curve from $\varphi(M)$ to $p$ since $\varphi(M)$ is closed and $N$ complete; this is a normal geodesic. This implies that $\eta$ is surjective.\par
In contrast to the cut locus of a point the cut distance function is in general not continuous. It  is upper semi-continuous and it is discontinuous, but not necessarily lower semi-continuous. In analogy to the cut locus of a point we have the following result:
\begin{prop}\label{focal-cut-point}
The cut locus only consists of focal and cut points.
\end{prop}
\prf We consider a $v\in \nu^1 M$ with $\sigma(v)<\infty$ such that $\sigma(v)v$ is not a
focal normal. We have to show that $\eta(\sigma(v)v)$ is a cut point. We construct
sequences $(\tilde{v}_n)$ in $\nu^1 M$ and $t_n>0$ such that
$t_n\tilde v_n$ is minimal with $\eta(t_n\tilde v_n)=\eta((\sigma(v)+1/n)v)$.
As $\eta|\bar B_r(\nu M)$ is proper for all $r\geq 0$ we can
assume that $t_n\tilde{v}_n$ converges to, say $t_0\tilde{v}$, where $\tilde v\in\nu^1 M$. Then
$\sigma(v)v$ and $t_0\tilde v$ are minimal and have the same endpoint. This implies $t_0=\sigma(v)$.
Since $\eta$ is injective on a neighborhood of $\sigma(v)v$ we have $\iota(\tilde{v})\neq \iota(v)$ and
$\eta(\sigma(v)\tilde{v})=\eta(\sigma(v)v)$.\eop

Later we will frequently us the following notion. If $r=\inf\{\sigma(v)\ |\ v\in\nu^1 M\}>0$ we call $r$ {\it injectivity radius} of $\varphi$ and $T_s=tube(M,s)=\exp B_s(\nu M)$ an {\it injectivity tube} of $M$ with radius $s$ for any $s$ with $0<s\leq r$. By definition for each $p\in T$ there is exactly one minimal normal $v$ with endpoint $p$ up to foot point. If $\varphi$ is injective, $\eta:B_s(\nu M)\to T_s$ is a diffeomorphism. With the aid of the cut locus we can now already give a necessary condition for a properly immersed submanifold $M$ with the two minimal conditions to give rise to a global foliation of the ambient space $N$. 
\begin{prop}\label{necessary_foliation}
Suppose a (topologically) closed submanifold $M$ satisfying the minimal conditions given in the introduction induces a global foliation $\F$. Then the cut distance function is constant along parallel normal fields.
\end{prop}
\prf Assume that $M$ induces a global foliation and that the cut distance function $\sigma$ is not constant for some parallel normal field. Then there is a $v_0\in\nu M$ with $\|v_0\|>\sigma(v_0/\|v_0\|)$ and a parallel translation $v_1\in\nu M$ with $\|v_1\|<\sigma(v_1/\|v_1\|)$, i.e. $v$ is minimal. We find a minimal vector $w_0\in \nu M$ with $\eta(w_0)=\eta(v_0)$. Since $M_v=M_w$ there is a normal parallel translation $w_1$ of $w_0$ with $\eta(w_1)=\eta(w_2)$. So $\|w_0\|<\|v_0\|=\|v_1\|\leq \|w_1\|$ in contradiction with $\|w_0\|=\|w_1\|$.\eop

Indeed the condition of the proposition is not fulfilled in the second example in \ref{bspfol}. This condition is general not sufficient. The search for a necessary and sufficient condition can be answered for the class of submanifolds with sections. We will do this in the next section.

\section{Submanifolds with Parallel Focal Structure}\label{sectionpar}
\subsection{General Properties} Let $M$ and $N$ be complete and connected Riemannian manifolds and $\varphi:M\to N$ an isometric immersion. The aim of this section is to find minimal conditions for $M$ to foliate the ambient space $N$ with parallel submanifolds.\par

For any $v\in T_x N$ we have a decomposition of $T_vTN$
$$
T_vTN=H^M_v\oplus V^N_v\cong T_xN\times T_xN
$$ 
into the horizontal space $H^N_v$ and the vertical space $V_v^N$ with respect to the Levi-Civita connection $\nabla$. The pullback of the metric on $N$ with the canonical isomorphism $H^M_v\oplus V^N_v\cong T_xN\times T_xN$ is the Sasaki-metric by definition. Let $M$ be a submanifold of $N$. Similar as above, for any $v\in\nu M$ we obtain a decomposition
$$
T_v\nu M=H^M_v\oplus V^M_v\cong T_xM\oplus\nu_x M
$$
into the horizontal space $H^M_v$ and the vertical space $V^M_v$ with respect to the normal Levi-Civita connection $\nabla^\perp$.
Obviously $V^M\subset V^N$ but in general we do not have $H^M\subset H^N$. Indeed, an element $\xi=(\xi_h,\xi_v)\in T_v\nu M$ is equal to $(\xi_h,\xi_v-A_v\xi_h)$ as an element of $T_vTN$, where $A$ is the shape operator of $M$.\par
We have the isomorphism between $T_vTN$ and the vector space of Jacobi fields of $N$ along 
$\gamma_v$ mapping an 
element $\xi\in T_vTN$ to the Jacobi field $J=J_\xi$ given by $(J(t),J'(t))=\phi^t_*(\xi_h,\xi_v)$, where $\phi^t$ is the time $t$ map of the geodesic flow $\phi:\RR\times TN\to TN$.
The inverse map is given by $J\mapsto (J(0),J'(0))$. The restriction of the first map 
to $T_v\nu M$ is an isomorphism onto the vector space $\CJ_M(v)$ of $M$-Jacobi fields 
along $\gamma_v$ 
$$
T_v\nu M\to\CJ_M(v);\quad \xi\mapsto J_\xi \quad\mbox{with}\quad
(J_\xi(t),J_\xi'(t))=\phi^t_*(\xi_h,\xi_v-A_v\xi_h).
$$
The inverse map is given by $J\mapsto 
(J(0),J'(0)^\perp)$. The decomposition $T_v\nu M=H^M_v \oplus V^M_v$ carries over to the 
decomposition of $\CJ_M(v)$ into a horizontal and a vertical subspace. We can describe a 
vertical/horizontal $M$-Jacobi field $J$ with initial condition $\xi\in T_v\nu M$ 
by a variational vector field. Define $V(s,t)=\eta(tX(s))$, where $X$ is a vector field along 
the constant curve $c\equiv x$ with $\frac{dX}{dt}|_{t=0}=\xi_v$ if $J$ is vertical, and 
a parallel normal field along $c$ in $M$ with $\dot c(0)=\xi_h$ if $J$ is 
horizontal. Then $J(t)=\partial_s V(0,t)$.\smallskip 

\begin{dfn}\label{sections} Let $\varphi:M\to N$ be an immersion, $\iota:\nu M\to TN$ be the canonical inclusion. For $x\in M$ we call an isometric immersion $i_x:\Sigma_x\to N$ (or shorter $\Sigma_x$) with $(i_x)_*(T_x\Sigma_x)=\iota(\nu_xM)$ a {\it section}, if it is totally geodesic in $N$ and if $\Sigma_x$ is complete. Note that if we compose $i_x$ with a covering onto $\Sigma_x$, we obtain another section. Since we want $i_x$ to be unambiguous, we also demand $y=z$ whenever $(i_x)_*(T_y\Sigma_x)=(i_x)_*(T_z\Sigma_x)$. The immersion $\varphi:M\to N$ is said to {\it admit sections} if $\Sigma_x$ is a section for every $x\in M$ and if there is exactly one section of $\varphi$ through every regular point of the normal exponential map, i.e. if $p\in i_x(\Sigma_x)\cap i_y(\Sigma_y)$ is regular then $i_x= i_y\circ \alpha$ for some isometry $\alpha:\Sigma_x\to\Sigma_y$.
\end{dfn}

In order to avoid a cumbersome notation, we use $\Sigma_x$ and the term section in two different ways. When it comes to point sets, for instance, if we write $p\in\Sigma_x$, we actually mean by $\Sigma_x$ the image of the immersion $i_x$. If we talk about tangent vectors or curves of $\Sigma_x$, i.e., if the context is a topological or differentiable one, we are of course referring to the underlying manifold structure of the section. This distinction is particularly important here, since we allow $\Sigma_x$ to have self-intersections.
\begin{lemma}\label{jacobifield}
Let $\gamma$ be a geodesic in a section $\Sigma=\Sigma_x$ with $\gamma(0)=p=\varphi(x)$. 
Then any Ja\-co\-bi field in $N$ along $\gamma$ can be decomposed into $J=J_1+J_2$, 
where $J_1$ is a Jacobi field of $\Sigma$ and $J_2$ is Jacobi field with $J_2(t)\in T_{\gamma(t)}\Sigma^\perp$ 
for every $t$. For an $M$-Jacobi field $J$ this decomposition is exactly the one into vertical 
and horizontal $M$-Jacobi fields along $\gamma$. In particular we have 
$J(t)\perp T_{\gamma(t)}\Sigma$ for a horizontal $M$-Jacobi field $J$.
\end{lemma}
\prf We write $J_1$ for the $T\Sigma$-part of $J$ and $J_2$ for the orthogonal part. 
Since $\Sigma$ is totally geodesic, the curvature operator $R_{\dot \gamma(t)}$ leaves 
$T_{\gamma(t)}\Sigma$ invariant and therefore, as a self-adjoint operator, also the orthogonal 
complement $T_{\gamma(t)}\Sigma^\perp$, so $R_{\dot\gamma(t)}J_1(t)\in T_{\gamma(t)}\Sigma$ and 
$R_{\dot\gamma(t)}J_2(t)\in T_{\gamma(t)}\Sigma^\perp$ . On the other hand we have 
$J_1''(t)\in T_{\gamma(t)}\Sigma$, since $\Sigma$ is totally geodesic, and 
$J_2''(t)\in T_{\gamma(t)}\Sigma^\perp$ for all $t$, because of 
$0=\frac{d^2}{dt^2}g(J_2(t),X(t))=g(J_2''(t),X(t))$ for any parallel field $X$ of $\Sigma$ along 
$\gamma$.
The Jacobi identity for $J$ gives 
$$
0=R_{\dot\gamma(t)}J(t)+J''(t)=(R_{\dot\gamma(t)}J_1(t)+J_1''(t))+(R_{\dot\gamma(t)}J_2(t)+J_2''(t)).
$$ 
Since the term in the first bracket lies in $T_{\gamma(t)}\Sigma$ and the term in the second in $T_{\gamma(t)}\Sigma^\perp$, the vector fields $J_1$ and $J_2$ are also Jacobi fields.
The second statement follows from the initial conditions $(J_i(0),J'_i(0)^\perp)$ of $J_i$ 
for $i=1,2$.\eop
 
The kernel of $d\eta(v)$ consists of $(J(0),J'(0)^\perp)$, where $J$ is an $M$-Jacobi field along 
$\gamma_v$ with $J(1)=0$. The decomposition $J=J_1+J_2$ as in the lemma
then implies that $\ker d\eta(v)$ is a direct                              
sum of a horizontal and a vertical subspace of $T_v\nu M$ and that the kernel of $d\eta(v)$ only has a non-trivial vertical component if and only 
if $\eta(v)$ is a conjugate point of $x$ along $\gamma_v$ in $\Sigma_{x}$.
Summing up, the decomposition of an $M$-Jacobi field $J$ into $J=J_1+J_2$ means that
\begin{eqnarray}
d\exp^\perp(v): H_v^M\oplus V_v^M&\to& T_{\eta(v)}\Sigma^\perp \oplus T_{\eta(v)}\Sigma\\
(J(0),J'(0)^\perp)&\mapsto&J_1(1)+J_2(1)\nonumber
\end{eqnarray}
splits as an orthogonal direct sum of linear maps $H_v^M\to T_{\eta(v)}\Sigma^\perp$ and
$V_v^M\to T_{\eta(v)}\Sigma$. We call this {\it splitting of $\eta$}.\par 
\begin{dfn} We call a focal normal $v$ of {\it horizontal/vertical type} if $\ker d\eta(v)$ has a 
non-trivial horizontal/vertical component. If a normal vector $v$ is not a focal normal of 
horizontal type we call $v$ {\it $f$-regular}. A point $p\in N$ is called {\it $f$-regular} 
if there is an $f$-regular normal $v$ such that $\eta(v)=p$. For a normal vector $v\in\nu_x M$
we call the dimension of the horizontal factor of $\ker d\eta(v)$ the {\it horizontal multiplicity} of 
$v$. Note that we have slightly changed the definitions given in \cite{ewert}.
\end{dfn}
\begin{dfn}\label{defsub}
An immersion $\varphi:M\to N$ has {\it parallel focal structure}, if
\begin{enumerate}
\item $\nu M$ is flat,
\item $\dim(\ker d\eta(v)\cap H^M_v)=\dim\ker d(\eta\circ v)$ 
is constant for any local parallel normal field $v$, i.e.
the horizontal focal data is invariant under normal parallel translation, and
\item $\varphi$ admits sections.
\end{enumerate}
\end{dfn}
In contrast to \cite{ewert} we do not demand the invariance of the vertical data. We will show in \Verweisp\ref{parallel-cutlocus} that this second invariance is an implication.\par
\begin{ex} 
Regular orbits of polar actions have parallel focal structure.  Isoparametric submanifolds in $\RR^{n+k}$ and equifocal submanifolds in simply connected, compact symmetric spaces obviously fulfill conditions (1) and (2) of a submanifold. The existence of sections for both classes of submanifolds follows from the properties of the respective ambient space. \Verweist\ref{maintheorem1} will show, that they admit sections if and only if the set of parallel manifolds builds a foliation on the regular set, which is known for both classes.
\end{ex}

We assume that $\varphi$ admits sections and that $\nu M$ is flat. We define two distributions $\CD$ and $\CD^\perp$ on the set $N_r$ of $f$-regular points in $N$ by $\CD^\perp(p)=T_p\Sigma$, where $\Sigma$ is a section through $p$; let 
$\CD$ be the orthogonal distribution. The distribution $\CD^\perp$ and therefore $\CD$ are well-defined on the set of regular points, since $M$ admits sections, but a priori not on the set of $f$-regular points. It is easy to see that both distributions are integrable on the regular set: Let $p$ be a regular point and $v\in \nu M$ with $\eta(v)=p$. Recall that $\nu M$ carries the horizontal foliation $\CP$ given by normal parallelity, and the vertical foliation $\CP^\perp$ given by the fibers of the projection $\nu M\to M$.  Now let $U_v$ be an open neighborhood of $v\in\nu M$ such that $\eta|U_v:U_v\to V$ from $U_v$ onto its image $V$ is a diffeomorphism. The map $\eta|U_v$ maps vertical leaves diffeomorphically onto open subsets of sections. The splitting of $\eta$ says that $d\eta$ maps the horizontal distribution on $\nu M$ to $\CD$, i.e. $d\eta(v)(T_xM)=\CD(\eta(v))$. Since $U_v$ is bifoliated and $\eta|U_v$ is a diffeomorphism, $V$ is also bifoliated with respect to $\CD$ and $\CD^\perp$.
We want to show that both distributions are also differentiable and well-defined on the set $N_r$ of $f$-regular points in $N$. 
Integrability is clear.
\begin{lemma}\label{uniquesection}
There is exactly one section $\Sigma$ through a given $f$-regular point $p$ 
and $\eta^{-1}(p)$ only consists of $f$-regular vectors that are tangential to $\Sigma$. Moreover, $N_r$ is open and dense in $N$ and there is a unique differentiable extension of $\CD^\perp$ on $N_r$. The distributions $\CD$ and $\CD^\perp$ give rise to a bifoliation $(\F_r,\F_r^\perp)$ of $N_r$.
\end{lemma}
\prf Existence follows by surjectivity of $\eta$. We show uniqueness. 
Let $v_0\in\nu_x M$ be an $f$-regular vector with $\eta(v_0)=p$. 
Then there is a simply connected neighborhood $U$ of $x$ in $M$ such that 
$(\eta\circ v)|U:U\to P_{v_0}=\eta(U)$ is a diffeomorphism, where $v$ is a parallel normal field on 
$U$ with $v_x=v_0$. We define $T=tube(P_{v_0},\eps)=\{\exp(\xi)\ |\ \xi\in B_\eps(\nu P_{v_0})\}$. 
By shrinking $U$ we can assume that $T$ is an injectivity tube around $P_{v_0}$ for small $\eps>0$. 
Let $\rho:T\to P_{v_0}$ be the projection.
We have $T_{\eta(v_z)}\Sigma_z\perp T_{\eta(v_z)}P_{v_0}$ for every $z\in U$ by 
the splitting of $\eta$. Therefore {\it a slice of the tube $T$ through $\eta(v_z)\in P_{v_0}$ 
coincides with the component of $\Sigma_z\cap T$ containing $\eta(v_z)$.} We can therefore extend $\CD^\perp$ differentiably to $T$ 
as the kernel of the differential of the submersion $\rho$. Since $\CD^\perp$ is defined on the open and 
dense set of regular points of $N$, 
{\it this extension is the unique differentiable extension of $\CD^\perp$.}\par
Let $w_0\in\nu_y M$ be another $f$-regular vector with $\eta(w_0)\in T$. The same process as for $v_0$ gives us a simply connected neighborhood $U'$ of $y$, a parallel normal field $w$ extending $w_0$, $P_{w_0}$ and its tube $T'$ with the same properties. By eventually shrinking $U'$ and the radius of $T'$ we can assume $T'\subset T$. By the uniqueness of a differentiable extension of $\CD^\perp$ we conclude that the slices of $T'$ are equal to the slices of $T$ intersected with the open set $T'$. In particular, if $\eta(w_0)=p$ this implies that $v_0$ and $w_0$ are tangential to the same section $\Sigma_x=\Sigma_y$. Since $w_0$ is $f$-regular, $\eta\circ w$ has maximal rank on a neighborhood of $y$. We can assume this neighborhood to be $U'$. Then $P_{w_0}$ intersects the slices transversally, i.e. $\rho\circ\eta\circ w:U'\to P_{v_0}$ is a diffeomorphism onto its image $(*)$.\par
We have seen above that the $f$-regular vectors in $\eta^{-1}(p)$ are tangential to the same section. Now we are going to show that any $w_0\in\nu M$ with $\eta(w_0)\in T$ is $f$-regular. Then $T\cap \Sigma$ is an open neighborhood of $p$ in $\Sigma$ only containing $f$-regular points. This implies that the set of $f$-regular points is open in $N$ and that $\eta^{-1}(p)$ only consists of $f$-regular vectors. We remark that this even shows that the $f$-regular points in a section $\Sigma$ are open in $\Sigma$ (see \Verweisc \ref{opendense}).  Let $w_0\in\nu_y M$ with $\eta(w_0)\in T$ and $U'$ a neighborhood of $w_0$ in $\nu M$ such that $\eta(U')\subset T$. We want to show that $w_0$ is $f$-regular. We can locally define a parallel normal field $w$ extending $w_0$. Then there is a simply connected neighborhood $U$ of $y$ in $M$ and an $\eps>0$ such that the image of $U$ under $(1+t)w$ lies in $U'$ for all $t\in(0,\eps)$ and such that $w'_z$ is $f$-regular for every $z\in U$, where $w'=(1+\eps)w$.
The geodesic $\gamma_{w(z)}$ intersects $P_{(1+\eps)w_0}$, the image of $\eta\circ w'$, orthogonally in $\gamma_{w(z)}(1+\eps)$ for all $z\in U$ by the splitting of $\eta$ or the Gauss Lemma for the normal exponential map. Then the image of $\gamma_w|[1,1+\eps]$ lies in a slice of the tube $T$. Therefore
$\rho\circ\eta\circ w=\rho\circ\eta\circ w'$
on $U$. Since the right side is a diffeomorphism $(*)$ this implies that also $\eta\circ w$ has maximal rank, i.e. $w_0$ is $f$-regular.\par
For any $f$-regular point $p$ we obtain a neighborhood $T$ as above that is bifoliated with respect to $\CD$ and $\CD^\perp$.
\eop

The lemma says that the preimage of a focal point $\eta(v)$, where $v$ is a focal normal of horizontal type, only consists of focal normals of horizontal type.\par
By the Theorem of Sard the set of regular points of $\eta$ is open and dense in $N$. Obviously the intersection of the set of regular points with $\Sigma$ is open in $\Sigma$. It is a priori not clear that the set of regular points in $\Sigma$ is dense in $\Sigma$. That this is true says the following corollary of the last lemma.
\begin{corollary}\label{opendense}
The subset of $f$-regular points in a section $\Sigma$ is open and dense in $\Sigma$.
\end{corollary} 
\prf We have seen in the proof of \Verweisl \ref{uniquesection} that the subset of $f$-regular points in $\Sigma$ is open in $\Sigma$. For the proof of density, see \cite{toeben}.\eop

Now let $\varphi$ be submanifold with parallel focal structure. Then every horizontal leaf $L_v$ in $\nu M$ through an $f$-regular vector $v$ is contained in $\eta^{-1}(N_r)$. The map $\eta:L_v\to M_v$ is a covering. Obviously $M_v$ is open in a leaf. Using completeness of $M$ one can easily see that $M_v$ is a leaf. The first statement of the following proposition follows from \Verweisl \ref{uniquesection}.
\begin{prop}\label{bifoliation}
Let $\varphi$ be a submanifold with parallel focal structure. Then the leaves of $\F_r^\perp$ are the parallel manifolds and the leaves of $\F_r^\perp$ the components of the restrictions of the sections to $N_r$. $\F_r^\perp$ is a totally geodesic foliation and therefore $\F_r$ a Riemannian foliation. Moreover, every parallel manifold has a flat normal bundle.
\end{prop}
\prf The characterization of the leaves of $\F_r$ is given above, the one for leaves of $\F_r^\perp$ is clear. By definition the sections are totally-geodesic, therefore $\F_r^\perp$ is a totally geodesic foliation. Let us consider a bifoliation $(\F_1,\F_2)$ of a Riemannian manifold, where the two foliations are orthogonal to each other. It is well-known that $\F_1$ is a Riemannian foliation if and only if $\F_2$ is a totally geodesic foliation. Let us consider this kind of bifoliation. Let $L$ be a leaf of $\F_1$ and $v\in\nu L$ with footpoint $p$. We consider a plaque $P$ through $p$ of a neighborhood $U$ of $p$ that is simple (or foliated) with respect to $\F_1$. There is a vector field $U$ tangential to $\F_2$ and foliated with respect to $\F_1$ extending $v$. The restriction of this vector field to $P$ is a parallel normal field of $P$. Therefore any leaf of $\F_1$ and in particular our parallel manifolds have a flat normal bundle. For more details see \cite{toeben}.\eop

Let $\bar{M}$ be the normal holonomy principal bundle over $M$ equipped with the metric such that the projection $\bar M\to M$ becomes a Riemannian 
covering. Its normal bundle is globally flat and $\bar M\to M$ has the lowest degree among all coverings of $M$ with this property. Each normal vector $v$ of $M$ canonically defines a global parallel normal field on $\bar{M}$, denoted by $\bar v$. We will denote the normal exponential map of $\bar M$ also by $\eta$. \Verweisp \ref{bifoliation} implies that $\varphi$ and 
$\eta \circ \bar v$ for $f$-regular $v$ factorize through injective immersions, the first even through an injective isometric immersion. 
If $\varphi$ is proper it factorizes finitely over an embedding; if in addition $v$ is $f$-regular and has finite normal holonomy degree then $\eta\circ \bar v: \bar M\to N$ is also a proper immersion, since $\eta$ restricted to $\bar B_r(\nu M)=\{w\in\nu M\mid \|w\|\leq r\}$ is proper for any $r\geq 0$. 
So from now on we can assume that $\varphi$ is injective and the inclusion map of $M$ into $N$. Let us repeat the definition of parallel and focal manifolds.
\begin{dfn} Let $\varphi:M\to N$ have parallel focal structure. We call 
$\eta\circ \bar v:\bar M\to N$ a {\it focal submanifold} of $M$ if $v\in\nu M$ is a 
focal normal of horizontal type, a {\it parallel submanifold},
if $v$ is $f$-regular. In any case we denote the image by $M_v$.
\end{dfn}

Let $v$ be a focal normal of horizontal type. Since the map $\eta\circ\bar v$ has constant rank, the set of connected components of preimages of $\eta\circ\bar v$ defines a foliation (the {\it focal foliation}) by the rank theorem which gives us simple sets for this foliation. The leaf through $x$ is called the {\it focal leaf} $F_{\bar v_x}$ through $x$ associated to $v$ (or to $\bar v_x$). If $\CG$ denotes the focal foliation, $\bar M/\CG$ endowed with the quotient topology is not necessarily Hausdorff or second countable but carries a natural differentiable structure by \Verweist VIII of \cite{palais} for which the map $\eta\circ\bar v:\bar M\to N$ induces an immersion $M/\CG\to N$. Thus $M_v$ is the image of an immersion. 

\subsection{The Blow-Up \boldmath $(\hat N,\hat\F,\hat\F^\perp)$\unboldmath}\label{secblowup} Each parallel manifold has the same set of sections as $M$ and thus we have a splitting of its normal exponential map. By a similar argument as in \Verweisl\ref{uniquesection}, we can show that each parallel submanifold $M'$ has the same set of $f$-regular points in $N$, namely $N_r$, and therefore admits sections by \Verweisp \ref{bifoliation}. The restriction $\eta'|L$ of the normal exponential map $\eta'$ of a parallel manifold $M'$ to a horizontal leaf $L'$ of the flat bundle $\nu M'$ through an $f$-regular vector is a covering map onto a leaf of $\F$. In order to show that $M'$ has parallel focal structure it remains to show that $\eta'|L$ has constant rank, if $L'$ is a horizontal leaf in $\nu M'$ through a focal normal of horizontal type. We will see this in \Verweisp\ref{stability}.\par

Our main goal in this section is to show first that $\F=\{M_v\ |\ v\in\nu M\}$ is a global foliation and then a singular Riemannian foliation. In this subsection we will associate to $(N,\F)$ a certain foliated manifold $(\hat N,\hat\F)$. An analysis of this foliation will yield the results. Boualem defines this Riemannian foliation $\hat\F$ in \cite{boualem} from a singular Riemannian foliation $\F$. Thus we cannot use his construction. Instead we build up $\hat\F$ with the normal exponential map.\par 

For an $f$-regular point $x\in N$ let $\eta_x:\nu M_x\to N$ be the normal exponential map of the leaf $M_x$. We define
$$
\hat\eta_x:\nu M_x\to G_k(TN); v\mapsto T_{\eta_x(v)}\Sigma_{y},
$$
where $y$ is the footpoint of $v$. Note that $\Sigma$ can have self-intersections in focal points of horizontal type. Therefore we only have a well-defined tangential space $T_p\Sigma$ in $f$-regular points $p$, so our above definition of $\hat\eta_x$ is not precise. A correct definition is as follows. Let $i:\Sigma_y\to N$ be the unique section through $y$. We identify any $z\in \Sigma_y$ with its image $i(z)$, if $i(z)$ is $f$-regular. We define $\hat\eta_x(v):=i_*(T_{\gamma_{v}(1)}\Sigma_y)$ where $\gamma_{v}$ is the geodesic in $\Sigma_y$ with initial vector $v$. Since this notation is too cumbersome, we prefer the first one, but we have to keep in mind that the expression $T_{\eta(v)}\Sigma_y$ depends on $v$ and not only on $\eta(v)$. 
Let
$$
\hat N=\{T_q\Sigma\ |\ \Sigma\ \mbox{is a section}, q\in\Sigma\}.
$$
and let $\hat\pi:\hat N\to N$ be the footpoint map of $G_k(TN)$ restricted to $\hat N$. 
Then we have $\hat N=\hat\eta_x(\nu M_{x})$ for any $f$-regular point $x\in N$ since the set of sections of two different parallel manifolds coincide. Also note that $\eta_x=\hat\pi\circ\hat\eta_x$ Our next aim is to give a bifoliated manifold structure to $\hat N$. The idea is to model $\hat N$ on the normal bundles of the parallel submanifolds, the charts being the maps $\hat\eta_x$. The normal bundle $\nu M$ has two natural, complemetary foliations $\CP$ and $\CP^\perp$, one given by the flat horizontal structure, the other by the fibers of the projection $\nu M\to M$. \par

Let $p\in N$ be arbitrary. We fix $r>0$ and take $\eps'>0$ to be smaller than the injectivity radius of any point $q\in \bar B_r(p)$ in $N$. There is an $f$-regular point $x$ and a vector $v\in\nu_x M_x$ with $\eta_x(v)=p$ that is not a focal normal of vertical type. One can see that $d\hat\eta_x(w)|\CH^M_w$ is injective for any $w\in\nu M_x$. Therefore $\hat\eta_x$ has maximal rank on a neighborhood of $v$, even if $v$ is a focal normal of horizontal type. This means there is a neighborhood $U$ of $v$ in $\nu M_x$ such that $\hat\eta_x|U:U\to G_k(TN)$ is an embedding into $G_k(TN)$ and such that the footpoint set $V$ of $\hat V:=\hat\eta_x(U)$ is contained in $\bar B_{\eps'}(p)$. We take a ball neighborhood $P$ of $x$ in $M_x$ and a neighborhood $U_0$ of $v$ in $\nu_x M_x$ such that $\phi:P\times U_0\to U; (y,w)\to w_y$ is an injective immersion into $U$, where $w_y$ is the normal parallel displacement of $w$ to $y$. We reduce $U$ to the image of $\phi$ so that $\phi$ becomes a diffeomorphism onto $U$. We choose an $f$-regular point $p'$ in $B^{\Sigma_x}_{\eps'}(p)$, such that $p\in B^{\Sigma_x}_\eps(p')$ for some $\eps$ with $0<\eps<\eps'$. The map $\eta_x|\phi(\{y\}\times U_0)$ is a diffeomorphism onto its image $V_y$ for any $y\in L$ by choice of $U$ (note that $U$ does not contain any focal normals of vertical type). We shrink $U_0$ such that this map is a diffeomorphism onto $V_x=B_\eps^{\Sigma_x}(p')$ for $y=x$. 

\begin{lemma}\label{localisometry}
The map $\alpha_y:V_x\to V_y;\eta_x(v_x)\mapsto \eta_x(v_y)$ is an isometry, where $v_x\in U_0$ and $v_y$ is the normal parallel displacement of $v_x$ to $y\in P$. 
\end{lemma}
\prf The set $V_r=V\cap N_r$ is open and dense in $V$ and $U_r=\eta_x^{-1}(V_r)$, saturated by leaves of the shape $P\times \{w\}, w\in U_0$, is open and dense in $U$. We consider the diffeomorphism $\eta:U_r\to V_r$. The bifoliation on $U_r$ is mapped to the bifoliation $(\F_r,\F_r^\perp)$ restricted to $V_r$. A normal foliated field on $U_r$ maps to a normal foliated field on $V_r$. This is a parallel normal field when restricted to the plaques of $\F_r|V_r$. Moreover, any such parallel normal field along a regular plaque is given this way. If $w\in U_0$ is  $f$-regular, $P_w=\eta(\phi(P\times\{w\}))$ and $X$ is a parallel normal field on $P_w$, then $\|(\alpha_y)_*X(\eta(w))\|=\|X(\alpha_y(\eta(w)))\|=\|X(\eta(w))\|$. It follows that $\alpha_y:V_x\cap N_r\to V_y$ is a local isometry. As $V_x\cap N_r$ is open and dense in $V_x$, $\alpha_y:V_x\to V_y$ is an isometry.\eop

There is exactly one $v'\in U_0$ with $\eta(v')=p'$. Let $P':=\eta_x(\phi(P\times \{v'\}))$. We define the diffeomorphism $h:P\to P'; y\mapsto \eta_{x}(\phi(y,v'))$. Similarly as for $U$ we have a natural diffeomorphism $\phi':P'\times U_0'\to B_\eps(\nu P')$. By the splitting of $\eta_x$, $\hat\eta_x(\phi(y,v'))=\nu_{h(y)}P'$. The map $\eta_{p'}:B_\eps(\nu_q P')\to B_\eps^{\Sigma_q}(q)$ is a diffeomorphism for any $q\in P'$ by the choice of $\eps$ and $B_\eps^{\Sigma_{h(y)}}(h(y))=V_y$ for any $y\in P$. Then $\hat\eta_{p'}\circ (\phi'(\{h(y)\}\times U_0'))$ is equal to the transversal plaque $\hat V_y$ for any $y\in P$ $(*)$. Moreover, the map $k:U_0\to U_0'$ defined by $k(w)=(\eta_{p'}|\phi'(\{p'\}\times U_0'))^{-1}(\eta_x(\phi(x,w)))$ is a diffeomorphism.
Now let $w\in U_0$ be an arbitrary $f$-regular vector and $u=k(w)\in U_0'$. We extend $w$ and $u$ to parallel normal fields on $P$ respectively $P'$. The images of $\eta_x\circ w$ and $\eta_{p'}\circ u$ lie in the same plaque of $\F_r|V_r$ by the splitting of $\eta_x$ and $\eta_{p'}$. As $\hat\pi$ is injective on $\hat\pi^{-1}(N_r)$, the image of $\hat\eta_{p'}\circ u$ lies in the plaque $\hat\eta_x(\phi(P\times\{w\}))$ in $\hat V$.  Together with $(*)$ we have $\hat\eta\circ w=\hat\eta_{p'}\circ u\circ h$
on $P$. By continuity we have 
$$
\hat\eta_x\circ\phi(y,w)=\hat\eta_{p'}\circ\phi(h(y),k(w))
$$
for any $y\in P$ and $w\in U_0$.\par
So far we have the following. Given any $k$-plane $\xi\in \hat N$, any normal vector $v$ of a parallel manifold $M_x$ (where $x$ is the footpoint of $v$) with $\hat\eta_x(v)=\xi$, that is not a focal normal of vertical type, defines as above a neighborhood $\hat V$ of $\xi$. A chart is given by $\hat\eta_x:U\to \hat V$. The discussion above implies that any two chart domains $V$ intersect in open subsets of each other. So the union of topologies on the various neighborhoods $V$ forms a basis for the topology on $\hat N$, and $\hat N$ is a topological manifold. In addition we see that the change of coordinates $(h,k)$ is differentiable, so $\hat N$ carries a differentiable structure. Since $\hat\eta_x$ is also differentiable as a map into $G_k(TN)$, the differentiable structure is the unique one for which the inclusion $\hat N\to G_k(TN)$ is an immersion.
Moreover, the chart $\hat\eta_{x}:U\to \hat V$ induces two foliations on $\hat V$ that are complementary to each other. The leaves of the first are given by $\hat\eta_x(\phi(P\times\{*\}))$, the second by $\hat\eta_x(\phi(\{*\}\times U_0))$. A look at the change of coordinates $(h,k)$ reveals that these local foliations coincide on intersections. This gives us a (vertical) foliation $\hat\F$ and a complementary (horizontal) foliation $\hat\F^\perp$ on $\hat N$. We will state this result in the next proposition.\par
Since we have not yet defined a metric on $\hat N$, the denotation of $\hat\F^\perp$ has to be justified. The Grassmann bundle carries a canonical metric (see appendix in \cite{toeben}) for which the projection $G_k(TN)\to N$ is a Riemannian submersion, and the horizontal distribution of this bundle is given as follows. Let $\xi\in G_k(TN)$ be a $k$-plane through a point $p\in N$ spanned by an orthonormal $k$-frame $(v_1,\ldots,v_k)$. Then the horizontal lift $\tilde c$ of a curve $c$ in $N$ with $c(0)=p$ to $\xi$ is given by 
$$
\tilde c(t)=\mbox{span}\Big\{
({\overset{t}{\underset{0}\parallel}}c) v_1,\ldots, ({\overset{t}{\underset{0}\parallel}}c) v_k\Big\}. 
$$
In particular, the tangent bundle $T\Sigma$ of a totally geodesic submanifold $\Sigma$ of $N$ is horizontal with respect to the projection $G_k(TN)\to N$. We denote the pullback of this metric under $\iota$ by $\hat g$. Note that $\hat\pi|T\Sigma:T\Sigma\to\Sigma$ is then an isometry.
\begin{prop}\label{liftbifoliation} 
$\hat N$ carries a natural differentiable structure, for which the inclusion into $G_k(TN)$ is an immersion. Moreover $\hat N$ has a natural Rie\-mannian/totally geodesic bifoliation $(\hat\F,\hat\F^\perp)$ with respect to the pullback metric $\hat g$ of $\hat N$ in $G_k(TN)$. We have 
$$
\hat\F^\perp=\{T\Sigma\ |\ \Sigma\mbox{\ is a section of}\ M\}.
$$ 
\end{prop}
This proposition is a strengthening of Boualem's result in \cite{boualem}. He states it for some differentiable structure and some metric on $\hat N$. We prove it for the natural differentiable structure and metric $\hat g$. Moreover we do not need that the leaves are relatively compact. We will call $(\hat N,\hat\F,\hat\F^\perp)$ the {\it blow-up} of $(N,\F)$, when we have established that $\F$ is a singular Riemannian foliation.\medskip

\prf The statements about the differentiable structure and the existence of the bifoliation were derived above. The description of $\hat\F^\perp$ is clear. We only have to show orthogonality. Then, since the leaves of $\hat\F^\perp$ are totally geodesic, the duality implies that $\hat\F$ is a Riemannian foliation. We consider a chart $\hat\eta_{x}:U\to \hat V$. For $v\in U$ with footpoint $x$ and a horizontal vector $X\in T_vU$ and a vertical vector $Y\in T_vU$. We have
$$
\hat g\big(d\hat\eta(v)X,d\hat\eta(v)Y)\big)=g\big(d\eta(v)X,d\eta(v)Y)\big)=0.
$$
The first equality is valid because $d\hat\eta(v)Y\in T_{\hat\eta(v)}T\Sigma$ is horizontal for $\pi:G_k(TN)\to N$, $\hat\pi\circ\hat\eta=\eta$ and because $\pi$ is a Riemannian submersion. The second equality follows from $d\eta(v)Y\in T_{\eta(v)}\Sigma_x$ and $d\eta(v)X\perp T_{\eta(v)}\Sigma_x$ by the splitting of $\eta$. This implies that $\hat\F^\perp$ is the orthogonal foliation to $\hat\F$ with respect to $\hat g$.\eop

\begin{dfn} 
Let $\F_i$ be a partition of $N_i$ for $i=1,2$ into injectively immersed submanifolds. A map $f:(N_1,\F_1)\to (N_2,\F_2)$ is {\it foliated}, if it maps each element of $\F_1$ onto an element of $\F_2$.
\end{dfn}
From the discussion before \Verweisp \ref{liftbifoliation} we have the following:
\begin{lemma}\label{liftbifoliation2}
For an $f$-regular point $x$ the map $\hat\eta_x:\nu M_x\to \hat N$ is foliated with respect to the natural bifoliation $(\CP,\CP^\perp)$ on $\nu M_x$ and $(\hat\F,\hat\F^\perp)$ on $\hat N$.\eop
\end{lemma}
\begin{prop}[Stability]\label{stability}
If $M$ has parallel focal structure, so has every parallel manifold $M_v$. If $M$ is properly immersed (therefore embedded) and has finite normal holonomy, so does every parallel manifold.
\end{prop}
\prf Let $L$ be a leaf of $\hat\F$ in $\hat N$. We claim that $c(L):=\ker d(\hat\pi|L)(X)$ does not depend on the choice of $X\in L$. Let $x$ be $f$-regular and $v\in \nu\bar M_x$ with $\hat\eta(v)\in L$. Then the image of $\hat\eta_x\circ v$ is $L$, because $\hat\eta_x$ is foliated. We know that $\hat\eta_x\circ\bar v$ has maximal rank. We consider the formula $\eta_x\circ\bar v=(\hat\pi|L)\circ\hat\eta_x\circ\bar v$. For $x\in M$, the map $\eta_x\circ\bar v=\eta\circ\bar v$ has constant rank by assumption, so $c(L)$ is independent of $X\in L$ and equal to the horizontal multiplicity of $v$. Conversely, the formula now implies that $\eta_x\circ\bar v$ has constant rank for an arbitrary $f$-regular point $x$ and $v\in\nu_xM_x$ with $\hat\eta_x(v)\in L$. It follows that $M_x$ has parallel focal structure.\par
Since $\eta|\bar B_r(\nu M)$ is proper for any $r\geq 0$ and $M$ has finite normal holonomy, each parallel manifold is properly immersed. Being leaves (\Verweisp\ref{bifoliation}) they are embedded.
It remains to show that every parallel manifold of $M$ has finite normal holo\-nomy. Let $M_x$ be the parallel submanifold through a point $x$ and $\Gamma_x$ be the normal holonomy group of $M_x$ in $x$, acting on $\nu_x M_x$. Since the parallel manifolds are closed and embedded, the orbits of $\Gamma_x$ are discrete and compact, thus finite. By linearity $\Gamma_x$ is finite.  \eop 

Ewert states this result in Proposition 2.9 in \cite{ewert}, but his proof is not correct. In the fourth last line of p.\,20 he writes that $V_*{\partial_t}(1,\cdot,t)$ is a parallel normal field along the focal submanifold through $V(1,0,t)$. This is not true. He refers to \Verweisp 2.4,  \cite{ewert}, which is not correct if $M_z$ is a focal submanifold; take $x:=z\circ c$ for instance.\par
We define $\hat M_x=\hat\pi^{-1}(M_x)$ for $f$-regular $x\in N$. Note that $\hat\pi:\hat\pi^{-1}(N_r)\to N_r$ is a foliated isomorphism. If we already knew that $\F=\{M_v\ |\ v\in\nu M\}$ is a global foliation we would have that $\hat\pi:(\hat N,\hat\F)\to (N,\F)$ is foliated. Later we prove that $\hat M_x$ is connected also for the endpoint $x$ of a focal normal of horizontal type. This will show that $\F$ is a global foliation.\par

\subsection{Global Foliation}\label{secglobalfoliation} We know that the parallel submanifolds are injectively immersed and orthogonal to the sections in each point of intersection. So far this is not clear for the focal submanifolds. This means that neither $T_pM_v$ nor $\nu_pM_v$ is defined for $p\in M_v\subset N$.  Let $v\in\nu_x \bar M$ be a focal normal of horizontal type and $p=\eta(v)$. Let $F=F_{\bar v_x}$ be the focal leaf associated to $v$ containing $x$. Define $F'_v=\bar v(F)$ and $W=(d(\eta\circ\bar v)(x)(T_x\bar M))^\perp$. Up to this point we have not assumed properness of $\varphi$. 
\begin{lemma}\label{focallemma}
Let $\varphi:M\to N$ be a proper immersion with parallel focal structure and finite normal 
holonomy. Let $v\in\nu_x\bar M$ 
be a focal normal of horizontal type, $p=\eta(v)$ and $x$ $f$-regular.
Then $\bigcup_{y\in F}\hat\eta_x(\bar v_y)=W$.
\end{lemma}
\prf First we prove the inclusion from left to right. The rank theorem states that we can write $\eta\circ\bar v:\bar M\to N$ locally in coordinates as $(x_1,\ldots,x_n)\mapsto (x_1,\ldots,x_{n-\mu(v)},0,\ldots,0)$, where $\mu(v)$ is the horizontal multiplicity of $v$. Since the focal leaf $F$ is compact, we find a neighborhood $U$ of $F$ that is saturated by focal leaves such that $\eta\circ\bar v:U\to P$ is a fibration onto its image $P\subset M_v$. In particular $(d(\eta\circ\bar v)(y)(T_y\bar M))^\perp=W$ for every $y\in F$. The splitting of $\eta$ implies $A_y:=\hat\eta(\bar v_y)\subset W$ for every $y\in F$.\par
Now let $w\in W$ be arbitrary. We look for a $y\in F$ such that $w\in A_y$. $F'=\bar v(F)$ is compact since $\eta|\bar B_r(\nu M)$ is proper for any $r\geq 0$ and $\varphi$ has finite normal holonomy. The time one map $\phi^1$ of the geodesic flow maps $F'$ 
diffeomorphically onto a compact submanifold $F^1$ of $W$. Therefore we find a shortest ray $\gamma$ in $W$ from $F^1$ to $w$. 
Then $\gamma$ is orthogonal to $F^1$ in some point $v':=\phi^1(\bar v_y)$, $y\in F$. As we will soon see $T_{v'}(A_y)=\nu_{v'}F^1$ in $W$ (we have $A_y\subset W$), which implies that $\gamma$ and therefore $w$ lies in $A_y$. We want to show $T_{v'}(A_y)=\nu_{v'}F^1$. First we prove $T_{v'}(A_y)\subset\nu_{v'}F^1$. We have
$$
T_{v'}F^1=\{(0,J_\xi'(1))\ |\ \xi\in T_{\bar v_y}F'\},
$$
where, because $T_{v'}F^1$ consists of elements $d\phi^1(\bar v_y)\xi=(J_\xi(1),J_\xi'(1))=(0,J_\xi'(1))$ for $\xi\in T_{\bar v_y} F'$. This implies $T_{v'}F^1\subset T_{v'}T_pN=V^N_{v'}$. Since $\xi$ is horizontal, $J_\xi'(t)$ is orthogonal to $T_{\gamma_{\bar v_y}(t)}\Sigma_y=A_y$ by the splitting of $\eta$. So $T_{v'}(A_y)\perp T_{v'}F^1$ also in $T_{v'}T_pN=V_{v'}^N$ for the Sasaki metric, hence $T_{v'}(A_y)\subset\nu_{v'}F^1$, where we consider $F^1$ as a submanifold of $W$. Since $\dim F^1=\dim F=\mu(v)$ and $\dim W=\mu(v)+k$, where $\mu(v)$ is the horizontal multiplicity of $v$ and $k$ the codimension of $M$, we have $T_{v'}(A_y)=\nu_{v'}F^1$ by equality of dimensions.\eop

Let $x$ be $f$-regular, $\eta=\eta_x, M=M_x$. 
\begin{lemma}\label{fibration}
Let $v\in\nu_x\bar M$ be a focal normal of horizontal type that is not of vertical type. Then there is a neighborhood $O$ of $x$ that is saturated by focal leaves of $v$, a relatively compact set $P$, $\eps>0$, a neighborhood $U_0$ of $\bar v_x\in\nu_x\bar M$ such that
\begin{enumerate}
\item $\eta\circ\bar v:O\to P$ is a surjective trivial fibration whose fibers are the focal leaves, i.e. the trivialization has the shape $O\cong F_v\times P$.
\item The map $\tilde\eta:O\times U_0\to T(P,\eps); (y,\bar w_x)\mapsto \eta(\bar w_y)$ is onto the tube $T(P,\eps)$.
\item $F_{\bar w_y}\subset F_{\bar v_y}$ for any $(y,\bar w_x)\in O\times U_0$. This means that the focal foliation given by $\eta\circ\bar w$ is finer than the focal foliation given by $\eta\circ\bar v$. 
\item Each section through a point $q\in T$ also contains the unique point $p'$ in $P$ that is in the same slice as $q$.
\item Let $p\in P$ and $S_p$ be the slice in $T$ through $p$. Then $S_q \subset S_p$ for any $q\in S_p$. 
\end{enumerate} 
\end{lemma}
\prf Any normal parallel translation $v'$ of $v$ with footpoint $y$ has the same multiplicity for $\exp^{\Sigma_y}$ as $v$ for $\exp^{\Sigma_x}$ as a consequence of \Verweisl\ref{localisometry}. Thus $v'$ is not a focal normal of vertical type. As before there is a neighborhood $O$ of the focal leaf $F$ saturated by focal leaves and a neighborhood $U_0$ of $v$ in $\nu_x\bar M$ such that $\eta\circ\bar v|O$ is a fibration onto its image $P$, which we can assume to be relatively compact, with typical fiber $F$ and such that $\hat\eta\circ\phi: O\times U_0\to \hat N$ is a diffeomorphism onto its image, where $\phi:O\times U_0\to \nu M;(y,w)\mapsto \bar w_y$. In particular $\hat\eta\circ\phi|(\{y\}\times U_0)$ is a chart for $\hat\F$. Let $\eps>0$ be smaller than the injectivity radius $i_N(q)$ in $N$ for all $q\in P$. Shrinking $U_0$ we can then assume that $\eta\circ\phi|(\{y\}\times U_0)$ is a diffeomorphism onto its image $V_y$ for any $y\in F$. We can also assume that $V_x$ is the image of the $\eps$-ball in $\hat\eta(\bar v_x)\subset T_pN$ around the origin under $\exp$ (if $\Sigma$ had no self-intersections we would write $V_x=B_\eps^{\Sigma_x}(p)$). By \Verweisl\ref{localisometry}, $\alpha_y:V_x\to V_y$ is an isometry. As $\eta\circ\phi(y,\,\cdot\,)=\alpha_y\circ\tilde\eta\circ\phi(x,\,\cdot\,)$ we have that $\eta\circ\phi(\{y\}\times U_0)=V_y$ is the image of the $\eps$-ball in $\hat\eta(\bar v_y)\subset T_pN$ around the origin under $\exp$. Then $\eta\circ\phi:O\times U_0\to T$ is surjective because the slice $S_q$ of $P$ in $T$ through $q\in P$ is equal to 
$$
S_q=\bigcup\{V_y\ |\ y\ \mbox{is in the focal leaf associated to\ }v\ \mbox{through\ }y\}
$$ 
for any $q\in P$ by \Verweisl\ref{focallemma}. \par
Let $y\in O$ and $u\in U_0$ be arbitrary. Let $F'=F_{\bar u_y}$ be the focal leaf associated to $u$ through $y$. Let $q=\eta(\bar u_y)$ and $p'=\eta(\bar v_y)\in P$. We want to show that $F'$ is contained in the focal leaf $F$ associated to $v$ through $y$. This is clear if $u$ is $f$-regular. We assume that $u$ is a focal normal of horizontal type. Obviously $V_y$ contains $p'$ and $q$. There is a vector $w\in \hat\eta(\bar v_y)\subset \nu_{p'}P$ of length smaller than $\eps$ with endpoint $q$. For $z\in U$ we define $w_z=d\alpha_z(p')w$, where $\alpha_z:V_y\to V_z$ as above but with central point $p'$ instead of $p$. The endpoint $\alpha_z(q)$ of $w_z$ is still in $T(P,\eps)$ because $\|w_z\|=\|w\|<\eps$ for all $z\in O$. For all $z\in F'\subset O$ we have $q=\eta(\bar u_z)=\alpha_z(q)$, thus $w_z=w$ since $w$ is unique among the vectors of $\nu P$ of length smaller than $\eps$ with endpoint $q\in T$. Therefore $\eta(\bar v_z)=\alpha_z(p')=p'$ for all $z\in F'$, so $F'\subset F_{\bar v_y}$. (In other words, the foliation of focal leaves given by $\eta\circ\bar u$ is finer than the foliation of focal leaves given by $\eta\circ\bar v$.) Therefore 
$$
S_q\subset\bigcup_{y\in F'}V_y\subset \bigcup_{y\in F}V_y=S_{p'}.
$$ 
\eop

We can deduce that the set $\hat\pi^{-1}(p)$ of tangential spaces in $p$ of sections through $p$ is equal to $\hat\eta(F'_v)$, where $\eta(v)=p$. To see this let $X\in\hat\pi^{-1}(p)$. Let $Y$ be the image of a small ball in $X$ around the origin, such that $Y\subset T(P,\eps)$. Assume that $X$ is not contained in $W=\bigcup_{y\in F}\hat\eta(\bar v_y)$. Then there is an $f$-regular point $z\in Y$ that lies in a slice $S_q$ for $q\in P\backslash\{p\}$. Since there is only one section through an $f$-regular point, $Y$ has to lie in $S_q$ and thus cannot contain $p$, contradiction. 
Now we will show in the following theorem that $M_v$ is embedded with the help of the blow-up.\par
Compare the first of following statements with the weaker result of Corollary 2.14 in \cite{ewert}. That
corollary is based on Lemma 2.13, \cite{ewert} which is not proved correctly (see the first sentence of the proof). 
\begin{thm}\label{maintheorem1}
If $\varphi:M\to N$ is a proper immersion with parallel focal structure and finite normal 
holonomy, then
$\F=\{M_v\ |\ v\in \nu M\}$ is a transnormal global foliation and the leaves of $\F$ are closed, embedded and orthogonal to each section they meet. Moreover, $\F$ is a singular Riemannian foliation admitting sections.
\end{thm}
{\prf} Assume $\eta(v)=\eta(w)=:p$ for $v,w\in\nu\bar M$. We have to show $M_v=M_w$. If $p$ is $f$-regular, 
then $v$ and $w$ are tangential to the same section by \Verweisl\ref{uniquesection}, so $\hat\eta(v)=\hat\eta(w)$. By \Verweisp\ref{liftbifoliation2} it follows $\hat\eta\circ\bar v(\bar M)=\hat\eta\circ\bar w(\bar M)$ and therefore $M_v=M_w$ because $\hat\pi\circ\hat\eta=\eta$. Now let $p$ be a focal point of horizontal type. It is possible that $\hat\eta(v)\neq \hat\eta(w)$. Both elements lie in $\hat\pi^{-1}(p)=\bigcup_{y\in F}\hat\eta(\bar v_y)$. Therefore there is a normal parallel translation $v'$ of $v$ with $\hat\eta(v')=\hat\eta(w)$. As above we conclude $M_v=M_{v'}=M_w$, so $\F$ is a global foliation.\par
We already know that the parallel submanifolds are closed and embedded by \Verweisp \ref{stability}. Now we consider $M_v$, where $v$ is a focal normal of horizontal type, and $p\in M_v$ arbitrary. Assume $\eta(v)=\eta(w)$ for $v,w\in\nu\bar M$ with footpoint $x,y$. Then by \Verweisl\ref{fibration} we find neighborhoods $O_1,O_2$ of $x,y$ that are saturated with focal leaves for $\eta\circ\bar v$ respectively $\eta\circ\bar w$ such that $\hat\eta\circ\bar v|O_1$ and $\eta\circ\bar w|O_2$ have the same image. Then $\eta\circ\bar v$ and $\eta\circ\bar w$ have the same image. As $\eta\circ\bar v:\bar M\to N$ is proper, $M_v$ is closed and embedded. Now $\nu M_v$ is well-defined for any $p\in M_v$. The previous discussion showed that $\nu_p M_v$ is the union of all $T_p\Sigma$, where $\Sigma$ is a section through $p$. This implies that also a focal submanifold intersects each section it meets orthogonally. This also implies that $\F$ is transnormal.\par
To prove that $\F$ is a singular Riemannian foliation, it remains to show that $\Xi(\F)$ acts transitively at a given point $p$. This is clear for $f$-regular $p\in N$, since the set of $f$-regular points $N_r$ is foliated by $\F_r$. Therefore we can assume that $p$ is a focal point of horizontal type of $M$ and $v\in\nu_x\bar M$ with $\eta(v)=p$. We assume that $v$ is not a focal normal of vertical type, otherwise we replace $M$ by a parallel manifold. Now we use the same objects as in \Verweisl \ref{fibration}.
We want to define a distribution $\CD'$ of dimension $\dim P$ on $T=T(P,\eps)$ such that $\CD'(q)\subset T_qM_q$. Let $q\in T$ be arbitrary.
Let $S_q$ be a slice of $M_q$ through $q$. Then there is a unique point $p'\in P$ such that the slice $S_{p'}$ of $T$ through $p'$ contains $q$. 
We define $\CD'(q)=T_qS_{p'}^\perp$. Since the distribution tangential to the
slices is differentiable so is $\CD'$. Since $S_q\subset S_{p'}$ by \Verweisl\ref{fibration}, we have 
$\CD'(q)\subset T_qM_q$. Thus, for any $p\in P$ and $X_0\in T_pP$ there is a vector field $X$ of $\CD'$ in $T$ extending $X_0$. If $f:N\to\RR$ is a bump function with support in $U$ and $f(p)=1$ then $fX\in \Xi(\F)$ with $(fX)_p=X_0$. Since $p$ and $X_0$ were arbitrary, $\Xi(\F)$ acts transitively on $T\F$.\eop

We can now exploit the theory of singular Riemannian foliations for submanifolds with parallel focal structure. Implications will be given in the next section. The converse was proven in \cite{alexandrino}. We give a different proof in the next section, see \Verweist \ref{marcos}.
\begin{rem}
Due to \Verweisl\ref{liftbifoliation2} $\hat\eta:(\nu M,\CP,\CP^\perp)\to (\hat N,\hat\F,\hat\F^\perp)$ is foliated. Since $\F$ is a global foliation, $\eta:(\nu M,\CP)\to (N,\F)$ is foliated, and then, because of $\eta=\hat\pi\circ\hat\eta$, also $\hat\pi:(\hat N,\hat\F)\to (N,\F)$.
Each parallel submanifold has the same focal submanifolds by the theorem.\par
A focal point of horizontal type of $M$ is also a focal point with the same horizontal multiplicity of any other parallel submanifold of $M$ and vice versa. For a focal normal $v$ of horizontal type $\eta\circ v:\bar M\to M_v$ is locally trivial fibration by \Verweisl \ref{fibration}, and so is the restriction of $\hat\pi$ to $\hat M_p$ for a focal point $p$ of horizontal type.
\end{rem}
The starting point of our work was the question, under which conditions a submanifold $M$ in $N$ with minimal conditions (1) and (2) stated in the introduction induces a global foliation $\F$ by parallel and focal submanifolds. For such a submanifold $M$ that has in addition sections $\Sigma_x$ for every $x\in M$ a necessary and sufficient condition to induce a global foliation is that $M$ admits sections. Sufficiency is provided by \Verweist \ref{maintheorem1}. Necessity is clear. Otherwise there is a regular point $p$ and two sections $\Sigma_1$ and $\Sigma_2$ with $T_p\Sigma_1\neq T_p\Sigma_2$. Then there are two parallel manifolds $M_{v_i}$ with $T_p{M_{v_i}}\perp T_p\Sigma_i\ (i=1,2)$, thus $T_p{M_{v_1}}\neq T_p{M_{v_2}}$, contradicting that $M$ induces a global foliation.\par

We call the elements of $\F$ {\it leaves}. A leaf is called {\it regular} if its dimension is 
maximal in $\F$, otherwise {\it singular}. A regular leaf with non-trivial normal holonomy is called {\it exceptional}.\par
A point in $N$ is $f$-regular if and only if it is contained in a regular leaf of $\F$. This justifies the denotation: the $f$ in $f$-regular stands for {\it foliation}.

\section{Singular Riemannian Foliations}\label{sectionSRF}
\subsection{Parallel Focal Structure of Regular Leaves} Let $\F$ be a singular Riemannian foliation admitting sections of a complete Riemannian manifold $N$. (For an introduction to singular Riemannian foliations, see section 6 in \cite{molino})
Let $M$ be a regular leaf and $\eta:\nu M\to N$ be its normal exponential map. We recall that a foliated vector field normal to the foliation on a simple neighborhood of a regular point (for $\F$) are parallel normal fields when restricted to the (regular) leaves. We can derive that $\nu M$ is flat. Therefore $\nu M$ is endowed with a natural foliation of horizontal leaves. The existence of sections implies the splitting of $\eta$. Therefore we can speak of $f$-regular points and vectors and of focal normals of horizontal/vertical type. From \cite{molino} we know that the stratum of regular points is open and dense in $N$. The following lemma is not difficult to prove. The second statement follows from the first with \Verweisc\ref{opendense}.
\begin{lemma}[\cite{toeben}]\label{equivalence_regular}
A point of $N$ is $\F$-regular if and only if it is $f$-regular. In particular the subset of $\F$-regular points in a section is open and dense.\hfill$\Box$
\end{lemma}
Hence a regular leaf $M$ admits sections in the sense of section \ref{sectionpar}.
\begin{prop}\label{expfoliated}
The map $\eta:\nu M\to N$ is foliated and the restriction of $\eta$ to a horizontal leaf in $\nu M$ has constant rank.
\end{prop}
\prf Let $v\in\nu M$ with endpoint $p$ and footpoint $x$. 
We define 
$$
Z=\{w\in L_v\ |\ \eta(w)\in M_p\},
$$
where $L_v$ is the horizontal leaf of $\nu M$ through $v$. We want to show that $Z$ is open and closed in $L_v$ and therefore equal to $L_v$ by connectivity. Let $w\in Z$ with footpoint $y$ and $q=\eta(w)$. Let $P_q$ be a relatively compact open neighborhood of $q$ in $M_q$ and let $T$ be an injectivity tube around $P_y$ (which is a distinguished neighborhood of $P_q$ in the sense of 6.2 in \cite{molino}). Since $\Xi(\F)$ acts transitively on $T\F$ we can assume that each plaque in $T$ intersects each slice of $T$ and always transversally. Thus the restriction of the projection $\rho: T\to P_q$ to an arbitrary plaque in $T$ is a surjective submersion. We choose a positive number $t<1$ such that $tv$ is $f$-regular and $\gamma_w|[t,1]$ lies in $T$. We see that $\gamma_w$ intersects $P_q$ orthogonally for $t=1$ since $\F$ is transnormal and $\gamma_w|[t,1]$ lies in the slice of $T$ through $q$. The leaf $M_{tv}$ is regular by the previous lemma. Let $L'=L_{(1-t)\phi^t(v)}$ be the horizontal leaf in $\nu M_{tv}$ containing  $(1-t)\phi^t(v)$. Observe that the map $\alpha:L_v\to L';\xi\mapsto (1-t)\phi^t(\xi)$ is a diffeomorphism and that $\eta_x|L_v=(\eta_{\eta(tv)}\circ\alpha)|L_v$. This means that we can replace $M$ by $M_{tv}$ for our considerations and assume that $\gamma_w|[0,1]$ is contained in a slice of $T$, so in particular the footpoint $y$ of $w$ lies in $T$. Let $P_y$ be the connected component of $M_y$ in $T$ containing $y$. We define the function $r:T\backslash P_{q}\to \RR$ measuring the distance to $P_{q}$ and let $X=-\grad r$ be the negative of the radial vector field. Then $w=\|w\| X_y$. Note that $X|P_y$ is a normal vector field of $P_y$. The flow of $X$ is a family of homotheties in $T$ centered at $P_q$ which respects the singular Riemannian foliation by the Homothety Lemma (see \Verweisl 6.2 in \cite{molino}). Due to \Verweisp 2.2 in \cite{molino} $X$ is a foliated vector field on a neighborhood of $P_y$ in $T$. Thus $X|P_y$ is a normal parallel field of $P_y$ and the image of $P_y$ under $X$ is an open subset of the horizontal leaf $L_v$ in $\nu M_x$ containing $w$. We want to show that $(\eta\circ (\|w\|X))|P_y=\rho|P_y$ which implies that $Z$ is open in $L_v$. But this follows from the observation that $\phi_X(t,z)=\gamma_{X_z}(t)$ for $t\in[0,\|w\|)$ and $z\in P_x$ where $\phi_X$ is the flow of $X$; note that $\|w\|$ is the distance of $P_y$ and $P_q$. We remark that this implies that $\eta|L_v$ has constant rank and its image is open in $M_p$. Now let $w\notin Z$ with footpoint $y$ and endpoint $q$. By assumption $q\notin M_p$. As above we show that an open neighborhood of $w$ in $L_v$ is mapped to $M_q$ which is disjoint to $M_p$ by definition of $\F$. Therefore the complement of $Z$ is also open. Thus $\eta(L_v)\subset M_p$.\par 
We will now show $\eta(L_v)=M_p$. We have seen above that $\eta(L_v)$ is open in $M_p$. It suffices to show that $\eta(L_v)$ is also closed in $M_p$. Let $q$ be an arbitrary point on the boundary of $\eta(L_v)$ in $M_p$. We have to show $q\in\eta(L_v)$. There is an injectivity tube $T$ of some open neighborhood $P_q$ of $q$ in $M_q$. As $\Xi(\F)$ acts transitively on $T\F$, we can assume that any plaque in $T$ meets any slice of $P_q$, and always transversally. Now there is a $w\in L_v$ such that $\eta(w)\in P_q$. As above we can assume that the footpoint $y$ of $w$ is contained in $T$. Then we define $X=-\grad r$ on $T\backslash P_q$ and we have $w=\|w\|X_y$. The endpoint of $\|w\|X_{y'}$ for $y'\in P_y$ is the unique point in the intersection of $P_q$ and the slice of $P_q$ containing $y'$. Since $P_{y'}$ meets any slice of $P_q$, in particular the slice through $q$, we have $q\in\eta(L_v)$ and $\eta(L_v)=M_p$. \eop

As a direct consequence of the lemma we obtain the following theorem of M. Alexandrino.
\begin{thm}[\cite{alexandrino}]\label{marcos}
A regular leaf of a singular Riemannian foliation admitting sections of a complete Riemannian manifold has parallel focal structure.\eop
\end{thm}
\begin{rem}\label{focalleaf}
Let $M$ be a regular leaf and $v\in\nu M$ with footpoint $x$ and singular endpoint $p$. We want to see that there is a compact submanifold $F\subset M$ of dimension equal to the horizontal multiplicity of $v$ to which we can extend $v$ to a parallel normal field such that image of $\eta\circ v$ is $p$. As above we can assume that $x$ is contained in an injectivity tube $T$ of a relatively compact singular plaque $P_p$ and that $\gamma_v:[0,1]\to N$ is the minimal normal geodesic from $x$ to $p$. Let $F$ be a connected component of the intersection of a regular leaf $M$ with the slice through $p$ containing $x$. As in the proof of \Verweisl \ref{expfoliated}, $v$ can be extended to a parallel normal field on a neighborhood of $F$ which coincides with the restriction of a radial field up to a scalar. The image of $\eta\circ v|F$ is $p$ and the corresponding focal leaf in $\bar M$ covers $F$. We call $F$ focal leaf in $M$.
\end{rem}
The following is a slice theorem for singular Riemannian foliations admitting sections.
\begin{thm}[\cite{alexandrino}]\label{slicetheorem} 
Let $\F$ be a singular Riemannian foliation admitting sections of a complete Riemannian manifold $N$. Let $p\in N$, $B_\eps(0_p)$ be the ball of $0_p$ in $\nu_pM_p$ for a small radius $\eps$ and $S_p=\exp^\perp(B_\eps(0_p))$. Then
the restriction $\F|S_p$ is a singular Riemannian foliation admitting sections that is isomorphic to the restriction of an isoparametric partition $\F'$ of $\RR^m$ to a ball neighborhood of the origin with the same codimension. This isomorphism is given by $\exp^\perp:B_\eps(0_p)\to S_p$, and it maps flat sections of $\F'$ to sections of $\F$ restricted to $S_p$.\eop
\end{thm}
An isoparametric partition is invariant under homotheties, i.e. maps $h_\lambda:\RR^m\to\RR^m; x\mapsto \lambda x$, where $\lambda\in\RR$. Therefore one can recover the isoparametric partition from its restriction to an open neighborhood of the origin. 
An isoparametric family of submanifolds of $\RR^m$ is given as the level sets of a transnormal map. Therefore the isoparametric family in $\RR^m$ and $\F|S_q$ are proper singular Riemannian foliations, i.e., its leaves are closed and embedded. 

\subsection{Transversal Holonomy}
Let $(N,\F)$ be as in the previous section. Then by Theorem \ref{marcos} and section \ref{sectionpar}
$$
\hat N:=\{T_p\Sigma\ |\ p\in N, \Sigma \mbox{\  is a section of\ }\F\mbox{\ through\ }p \}
$$
carries the unique differentiable structure for which the inclusion $\hat N\to G_k(TN)$ is an immersion (see \Verweisp\ref{liftbifoliation}). Moreover, $\hat N$, endowed with the pull-back metric, carries by a Riemannian/totally geodesic bifoliation $(\hat\F,\hat\F^\perp)$, where 
$$
\hat\F^\perp=\{T\Sigma\ |\ \Sigma\ \mbox{is a section of}\ \F\}.
$$ 
We call $\hat\F$ {\it vertical} and $\hat\F^\perp$ {\it horizontal} foliation. Since $\hat\eta$ is foliated by \Verweisl \ref{liftbifoliation2} and $\eta$ by \Verweisl \ref{expfoliated}, so is the footpoint map $\hat\pi:(\hat N,\hat\F)\to (N,\F)$ because of $\eta=\hat\pi\circ\hat\eta$. The same arguments for $\CP^\perp$ and $\hat\F^\perp$ show that $\hat\pi$ maps a horizontal leaf $T\Sigma$ to the section $\Sigma$. Here by $\Sigma$ we mean the submanifold and not its image, which can have self-intersections in singular points. When we introduced the metric $\hat g$ on $\hat N$ in subsection \ref{secblowup}, we noticed that $\hat\pi:T\Sigma\to\Sigma$ is an isometry.\par
The discussion shows that $\hat M=\hat\pi^{-1}(M)$ is a leaf of $\hat\F$ for regular $M$, and a union of leaves if $M$ is singular. Let $M$ be singular. Now we want to see that $\hat M$ is also connected. In subsection \ref{secglobalfoliation} we used properness and finite normal holonomy assumptions on the submanifold with parallel focal structure to prove this (see the paragraph before \Verweist \ref{maintheorem1}). If a partition is already given, in our case $\F$, then these conditions are not necessary. Let $p$ in $M$. By definition $\hat\pi^{-1}(p)$ is the set of tangential spaces in $p$ of sections through $p$. It suffices to show that this set is contained in one leaf of $\hat\F$. Let $S_p$ be a slice through $p$. The corresponding isoparametric partition of $\nu_pM_p$ given by \Verweist\ref{slicetheorem} has closed and embedded regular leaves with parallel focal structure and finite normal holonomy. Let $L$ be a regular leaf of this isoparametric partition. Now \Verweisp\ref{focallemma} describes the set of sections through $p$ as the image of a focal leaf associated to $v$ under $\hat\eta\circ\bar v$ for some $v\in\nu L$, so $\hat\pi^{-1}(p)$ is contained in one leaf. This means $\hat M_p:=\hat\pi^{-1}(M_p)$ is a leaf. Therefore
$$
\hat\F=\{\hat\pi^{-1}(M)\ |\ M\in\F\}.
$$
For a curve $\tau:[0,1]\to N$ in a regular leaf of $\F$ and a curve $\sigma:[0,1]\to N$ in a section, both starting in an $\F$-regular point, we define the {\it lifts} $\hat\tau(t):=T_{\tau(t)}\Sigma_{\tau(t)}$ and $\hat\sigma(t):=T_{\sigma(t)}\Sigma_{\sigma(0)}$. Obviously $\hat\pi\circ\hat\tau=\tau$ and $\hat\pi\circ\hat\sigma=\sigma$.
\begin{lemma}\label{ehresmann2}
Let $x_0$ be $\F$-regular, let $\tau:[0,1]\to N$ be a curve in $M_{x_0}$ and $\sigma:[0,1]\to N$ be a curve in $\Sigma_{x_0}$ with $\tau(0)=\sigma(0)=x_0$. Then there is a unique continuous map $H=H_{(\tau,\sigma)}:[0,1]\times[0,1]\to N$ with
\begin{enumerate}\item $H(\,\cdot\,,0)=\tau$,
\item $H(0,\,\cdot\,)=\sigma$,
\item $H(\,\cdot\,,t)$ is contained in a leaf of $\F$,
\item $H(s,\,\cdot\,)$ is contained in a section.
\end{enumerate}
\end{lemma}
\prf First assume that $\F$ is a (regular) Riemannian foliation admitting sections. The set of sections forms a totally geodesic foliation orthogonal to $\F$. In this case the statements are due to \Verweisc 2.7 of \cite{blumenthalhebda}, which is based on \Verweisl 2.6. Note that for the proof of the latter one can drop completeness of $N$ and assume completeness of the horizontal leaves, the sections, instead. In particular the statements are valid for $(\hat\F,\hat\F^\perp)$.\par
Now let $\F$ be a singular Riemannian foliation with sections. Existence follows by $H_{(\tau,\sigma)}:=\hat\pi\circ H_{(\hat\tau,\hat\sigma)}$, where $\hat H_{(\hat\tau,\hat\sigma)}$ is defined as in the lemma for the bifoliation $(\hat\F,\hat\F^\perp)$ of $\hat N$. We want to show uniqueness of $H_{(\tau,\sigma)}$. Let $H$ be arbitrary with the four properties in the lemma. We define $\hat H(s,t):=T_{H(s,t)}\Sigma_{\tau(s)}$. Obviously $\hat H(s,\,\cdot\,)$ lies in the horizontal leaf $T\Sigma_{\tau(s)}$. The curve $H(\,\cdot\,,t)=\hat\pi\circ\hat H(\,\cdot\,,t)$ lies in the leaf $M_{\sigma(t)}$ by assumption. By the discussion at the beginning of this section, $\hat\pi^{-1}(M_{\sigma(t)})$ is a leaf of $\hat\F$. Therefore $\hat H(\,\cdot\,,t)$ is contained in a vertical leaf. By uniqueness we have $\hat H=\hat H_{(\hat\tau,\hat\sigma)}$ and therefore $H=\hat\pi\circ H_{(\hat\tau,\hat\sigma)}$ is also determined. \eop

 A continuous map $H:[0,1]\times [0,1]\to N$, such that $H(\,\cdot\,,t)$ is vertical for any $t$ and $H(s,\,\cdot\,)$ is horizontal for any $s$, is called {\it rectangle} with {\it initial vertical/horizontal curve} $H(\,\cdot\,,0)$/$H(0,\,\cdot\,)$, {\it terminal vertical/horizontal curve} $H(\,\cdot\,,1)$/$H(1,\,\cdot\,)$ and {\it diagonal} $t\mapsto H(t,t)$. We write $T_\sigma\tau=H_{(\tau,\sigma)}(\,\cdot\,,1)$ and $T_\tau\sigma=H_{(\tau,\sigma)}(1,\,\cdot\,)$. For a vertical curve $\tau$ in a leaf $M\in\F$ respectively a horizontal curve $\sigma$ in a section $\Sigma$ we write $[\tau]$ respectively $[\sigma]$ for the equivalence class of curves under homotopy in $M$ respectively $\Sigma$ fixing endpoints. Then $[T_\sigma\tau]$ and $[T_\tau\sigma]$ only depend on $[\tau]$ and $[\sigma]$. We remark that for any curve $\mu:[0,1]\to N$ we find a unique rectangle $H:[0,1]\times [0,1]\to N$ with diagonal $\mu$. We write $\mu_v$ respectively $\mu_h$ for the initial vertical respectively horizontal curve of $H$ and $\mu^v$ respectively $\mu^h$ for the terminal vertical respectively horizontal curve of $H$. In the sequel we make the following convention: We write $H_{(\tau,\sigma)}$ for rectangles in $N$ and, as in the proof of \Verweisl \ref{ehresmann2}, $\hat H_{(\hat\tau,\hat\sigma)}$ for rectangles in $\hat N$ with respect to the bifoliation $(\hat\F,\hat\F^\perp)$. Then $H_{(\tau,\sigma)}=\hat\pi\circ\hat H_{(\hat\tau,\hat\sigma)}$ as in the proof.\par
We recall that the universal cover $\tilde M$ of a manifold $M$ is equal to the set of equivalence classes of curves starting from a fixed point $x_0$, where the equivalence is given by homotopy fixing endpoints; in some cases we write more precisely $\widetilde{(M,x_0)}$. The covering map $\tilde M\to M$ is given by $[\sigma]\mapsto \sigma(1)$. Let $x_0$ be arbitrary and let $M\in\F$ the leaf through $x_0$ and $\Sigma$ the section through $x_0$. Then
\begin{eqnarray*}
\tilde M&=&\{[\tau]\ |\ \tau\ \mbox{is vertical and}\ \tau(0)=x_0\}\\
\tilde \Sigma&=&\{[\sigma]\ |\ \sigma\ \mbox{is horizontal and}\ \sigma(0)=x_0\}\\
\bar N&=&\{[\mu]\ |\ \mu\ \mbox{is a curve in}\ \hat N\ \mbox{and}\ \mu(0)=x_0\},
\end{eqnarray*}
where $\bar N$ denotes the universal covering of $\hat N$ and $x_0$ is identified with its lift $T_{x_0}\Sigma_{x_0}$. The manifold $\bar N$ is endowed with the pull-back bifoliation of the covering $\bar N\to \bar N$ and $\tilde M\times\tilde\Sigma$ carries a natural bifoliation. Due to \cite{blumenthalhebda2} the map $\Phi:\tilde M\times\tilde\Sigma\to \bar N; ([\tau],[\sigma])\mapsto [t\mapsto\hat H_{(\hat\tau,\hat\sigma)}(t,t)]$ (the diagonal) is a bifoliated diffeomorphism (i.e. foliated with respect to both pairs of foliations) with inverse map $[\mu]\mapsto ([\mu_v],[\mu_h])$ (see \cite{toeben} for details). Consequently
\begin{eqnarray*}
\Psi:\tilde M\times\tilde\Sigma&\to& \hat N\\
([\tau],[\sigma])&\mapsto&\hat H_{(\hat\tau,\hat\sigma)}(1,1)
\end{eqnarray*}
is a bifoliated universal covering map of $\hat N$. We define
\begin{eqnarray*}
\psi:\tilde M\times\tilde\Sigma&\to&N\\
([\tau],[\sigma])&\mapsto&H_{(\tau,\sigma)}(1,1).
\end{eqnarray*}
By definition 
$$
\hat\pi\circ\Psi=\psi.
$$
Since the footpoint map $\hat\pi:(\hat N,\hat\F)\to (N,\F)$ is foliated and maps a horizontal leaf $T\Sigma$ isometrically to the corresponding section $\Sigma$, it follows that $\psi$ is foliated and maps horizontal leaves onto sections. We sum up:
\begin{prop}\label{char2}
The map $\Psi$ is the universal covering map, and it is foliated with respect to the natural bifoliation of $\tilde M\times\tilde\Sigma$ and to $(\hat N;\hat\F,\hat\F^\perp)$. The map $\psi$ is foliated with respect to the vertical foliation on $\tilde M\times \tilde\Sigma$ and $(N,\F)$ and its restriction to a horizontal leaf is a Riemannian covering to a section.\eop
\end{prop}
\begin{rem} The map $\psi$ completely describes the singular Riemannian foliation $\F$ of $N$. The singular values of $\psi$ are exactly the singularities of $\F$. It is a covering when restricted to the regular set.
\end{rem}
Let $\tau:[0,1]\to M$ be a curve with $\tau(0)=x_0$. The map $T_{\hat\tau}:\widetilde{(T\Sigma,x_0)}\to \widetilde{(T\Sigma,\hat\tau(1))};[\sigma]\mapsto [T_{\hat\tau}{\hat\sigma}]$ defined  with respect to $(\hat\F,\hat\F^\perp)$ is an isometry due to \cite{blumenthalhebda2}. On the other hand we have $T_\tau:\widetilde{(\Sigma,x_0)}\to \widetilde{(\Sigma,\tau(1))};\ [\sigma]\mapsto [T_\tau\sigma]$ with respect to $\F$ and the family of sections. Since the isometry $\hat\pi:T\Sigma\to\Sigma$ passes to an isometry of the universal covers, we can identify $T_{\hat\tau}$ and $T_\tau$. In particular:
\begin{prop}
The sections have the same Riemannian universal cover. Similarly the regular leaves of $\F$ have the same universal cover.\eop
\end{prop}
This proposition describes a topological difference between a singular Riemannian foliation admitting sections and a polar action, namely the normal holonomy of a section. While the sections of a polar action are isometric to each other, the sections of a singular Riemannian foliation only have the same Riemannian universal cover. We want to explain this in more detail. We can define a local isometry along a vertical curve $\tau$ starting in $\Sigma$ similarly as in Verweisl \ref{localisometry}. It is important to know that in general such a map cannot be extended to an isometry that is defined on all of $\Sigma$. For instance consider the Klein bottle $N=[0,1]^2/\sim$, where we identify the two vertical edges in opposite direction and the horizontal ones in common direction. The two partitions, the one into vertical, the other into horizontal lines, build a Riemannian/totally geodesic bifoliation, so in particular a singular Riemannian foliation admitting sections. Take $M$ to be a vertical line and $\Sigma$ to be the exceptional horizontal line. Let $\tau$ be a curve in $M$ from a point in $\Sigma$ to a point that is not in $\Sigma$. Obviously we cannot extend a local isometry defined as above to a map defined on $\Sigma$ that respects the foliation. But we can develop these maps on the Riemannian universal cover $\tilde\Sigma$ and this is $T_{[\tau]}:\widetilde{(\Sigma,{\tau(0)})}\to \widetilde{(\Sigma,{\tau(1)})}$. The set of the above local isometries along vertical curves $\tau$ that start and end in $\Sigma$ becomes a pseudogroup of local isometries on $\Sigma$ while the set of its developments $T_{[\tau]}$ becomes a group acting on $\tilde\Sigma$ that we will later denote by $\tilde\Gamma$. It is more convenient to work with this group than with the corresponding pseudogroup. On the other hand we have to handle additional elements, namely the deck transformations of $\pi_\Sigma:\tilde\Sigma\to\Sigma$, which are contained in $\tilde\Gamma$, but do not contribute to the geometry of $\F$. For a special choice of sections we can divide them out and obtain a group $\Gamma$ acting on $\Sigma$, that completely describes the holonomy of $\F$. In the case that $\F$ is the orbit decomposition of a polar action, $\Gamma$ is the generalized Weyl group of a polar action.\par

For a vertical/horizontal  curve $c$ and a horizontal/vertical curve $d$ (horizontal means lying in a section) starting in the same regular point $x_0$ we denote by $T_cd$ the terminal horizontal/vertical edge of the homotopy $H_{(c,d)}$. 
$$
\tilde\Lambda=
\left\{[\tau]T_{[\sigma]}:\tilde M\to\tilde M \bigg|
\begin{array}{ll}
\ &\tau\ \mbox{is vertical}, \sigma\ \mbox{is horizontal}\\
&\tau(0)= \sigma(0)=x_0,\sigma(1)=\tau(1)
\end{array}
\right\}
$$
and
$$
\tilde\Gamma=
\left\{[\sigma]T_{[\tau]}:\tilde\Sigma\to\tilde\Sigma\ \bigg|
\begin{array}{ll}
\ &\tau\ \mbox{is vertical}, \sigma\ \mbox{is horizontal}\\
&\tau(0)= \sigma(0)=x_0,\sigma(1)=\tau(1)
\end{array}
\right\}.
$$
Later we focus on $\tilde\Gamma$. It turns out that this group (see below) carries the information about the tranversal geometry of $\F$. As we noted before, identifying $\Sigma$ with $T\Sigma$ and their universal covers, $T_\tau$ defined with respect to $\F$ and the family of sections and $T_{\hat\tau}$ defined with respect to $(\hat\F,\hat\F^\perp)$ are the same; here it is important that we have chosen $x_0$ to be regular (check the assumption in \Verweisl \ref{ehresmann2}).  As a consequence we can identify the above two sets with their counterparts for $(\hat\F,\hat\F^\perp)$. For $\tilde\Gamma$ this means that the transversal geometry of $\F$ can be read off from that of its blow-up $\hat\F$. The rules in the next lemma work for $T$ of $(\hat\F,\hat\F^\perp)$. For $T$ of $\F$ the second rule has an exception: If $c_2$ is horizontal and $c_2(0)$ is singular, the assumption in \Verweisl \ref{ehresmann2} are not met and $T_{(\cdot)}c_2$ is not defined. This case will not occur in our discussion.
\begin{lemma}\label{recrules} Let $c_1$ and $c_2$ both be horizontal/vertical with $c_1(0)=x_0$ and $c_1(1)=c_2(0)$ and $d$ vertical/hori\-zontal with $d(0)=c_1(0)$. Then
$$
T_{c_1c_2}=T_{c_2}\circ T_{c_1}\qquad \mbox{and}\qquad T_{d}(c_1c_2)=T_dc_1\cdot T_{T_{c_1}d}c_2.
$$
\end{lemma}
\prf The proof is clear.\eop

\begin{lemma}
$\tilde\Gamma$ is a subgroup of $I(\tilde\Sigma)$ and $\tilde\Lambda$ a subgroup of $\mbox{Diff}(\tilde M)$.
\end{lemma}
\prf This is an easy application of \Verweisl \ref{recrules}. See \cite{toeben}.\eop

We have already explained the geometric meaning of these groups. Now we want to give a relation between $\pi_1(\hat N,x_0), \tilde\Gamma$ and $\tilde\Lambda$. We identify the actions of $\tilde\Gamma$ and $\tilde\Lambda$ with the corresponding actions for $(\hat\F,\hat\F^\perp)$. For $[\mu]\in\pi_1(\hat N,x_0)$ we define $\tilde\gamma_{[\mu]}:= [(\mu_h)]T_{[(\mu^v)^{-1}]}\in \tilde\Gamma$ and $\tilde\lambda_{[\mu]}:= [(\mu_v)]T_{[(\mu^h)^{-1}]}\in \tilde\Lambda$, where $\mu_v$ is the initial vertical edge of the unique rectangle that has diagonal $\mu$, and so on. One can easily show that $\tilde\gamma:\pi_1(\hat N,x_0)\to\tilde\Gamma; [\mu]\mapsto \tilde\gamma_{[\mu]}$ and $\tilde\lambda:\pi_1(\hat N,x_0)\to\tilde\Lambda; [\mu]\mapsto \tilde\lambda_{[\mu]}$ are homomorphisms.
$\pi_1(\hat N,x_0)$ acts naturally from the left on the universal cover of $\hat N$. We want to transfer this action to $\tilde M\times\tilde \Sigma$ via the foliated isomorphism $\Phi$. Let $[\mu]\in \pi_1(\hat N,x_0)$ and $\Phi([\tau],[\sigma])=[\nu]$. Then 
$$
\Phi^{-1}([\mu][\nu])=([(\mu_h)]T_{[(\mu^v)^{-1}]}[\tau],[\mu_v]T_{[(\mu^h)^{-1}]}[\sigma])=(\tilde\gamma_{[\mu]}([\tau]),\tilde\lambda_{[\mu]}([\sigma])).
$$
This shows that the action of $\pi_1(\hat N,x_0)$ respects the natural bifoliation on $\tilde M\times\tilde\Sigma$. Thus: 
\begin{prop} 
$\pi_1(\hat N,x_0)$ is a subgroup of $\mbox{Diff}(\tilde M)\times \mbox{I}(\tilde\Sigma)$. The projection of $\pi_1(\hat N,x_0)$ on the first component is $\tilde\Gamma$, the one on the second is $\tilde\Lambda$. \eop
\end{prop}
The projection homomorphisms are $\tilde\gamma$ and $\tilde\lambda$. This describes a new view on the transversal holonomy group, even for the Weyl group of a polar action, for instance for the isotropy representation of a symmetric space. In \cite{toeben} we also show for a bifoliation $(\hat N,\hat\F,\hat\F^\perp)$: 
\begin{prop}
If $\hat M$ is a leaf of $\hat\F$ and $\Sigma$ a leaf of $\hat\F^\perp$, then \medskip

\hfill$|\pi_1(\hat N,x_0)|=|\pi_1(\hat M,x_0)|\cdot |\pi_1(\Sigma,x_0)|\cdot |\hat M\cap \Sigma|.
$\hfill\eop
\end{prop}
In the case of infinity, the interpretation of this equation is that the left value is infinity if and only if at least one of the factors on the right side is infinity.\par
Now we want to focus on $\tilde\Gamma$. There is a natural injective representation $\rho_\Sigma:\pi_1(\Sigma,x_0)\to \tilde\Gamma;[\sigma]\mapsto [\sigma]$. This means that $\tilde\Gamma$ contains the deck transformations of $\pi_\Sigma:\tilde\Sigma\to\Sigma;[\sigma]\mapsto \sigma(1)$. Since the elements of $\pi_1(\Sigma,x_0)$ do not contribute to the holonomy action on $\Sigma$ when restricting them to local transformations on $\Sigma$, it is natural to ask when we can divide $\pi_1(\Sigma,x_0)$ out of $\tilde\Gamma$. The following lemma gives a geometric condition.
\begin{lemma}
$\tilde\Gamma$ normalizes $\pi_1(\Sigma,x_0)$ if $\Sigma$ has trivial normal holonomy, i.e. if $T_g=\id_{\tilde M}$ for any $g\in\pi_1(\Sigma,x_0)$. 
\end{lemma}
\prf Let $\gamma=[\sigma]T_{\tau}\in\tilde\Gamma$ and $g\in\pi_1(\Sigma,x_0)$. Then for any $\lambda\in\tilde\Sigma$ we have 
$$\begin{array}{rl}
\gamma(g\lambda)=&[\sigma]T_{[\tau]}(g\lambda)\\
=&[\sigma]T_{[\tau]}g\cdot T_{T_{g}[\tau]}(\lambda)\\
=&[\sigma]T_{[\tau]}g[\sigma]^{-1}\cdot [\sigma]T_{[\tau]}(\lambda)\qquad \mbox{because}\ T_g=\id_{\tilde M}\\
=&g'\gamma(\lambda)
\end{array}$$
for $g'=[\sigma]T_{[\tau]}g[\sigma]^{-1}\in\pi_1(\Sigma,x_0)$.\eop

In the case of the lemma, we say that $\Sigma$ is an {\it exceptional section}. We remark that this is the generic case if each section is embedded. We call the group
$$
\Gamma=\tilde\Gamma/\pi_1(\Sigma,x_0)
$$
the {\it transversal holonomy group} of $\Sigma$. It is a subgroup of $I(\Sigma)$. The transversal holonomy group generalizes the Weyl group of the isotropy representation of a symmetric space (or polar action) and the fundamental domains of $\Gamma$ generalize the Weyl chambers. Note that $\Gamma$ is independent of the choice of the fixed regular point $x_0$ in the given section. But it depends on the choice of the section $\Sigma$, unlike $\tilde\Gamma$. Similarly we define the subgroup $\Lambda=\tilde\Lambda/\pi_1(M,x_0)$ in $\mbox{Diff}(\tilde M)$ if $M$ has trivial normal holonomy.\par
Moreover there is a representation $\rho_M:\pi_1(M,x_0)\to \tilde\Gamma;[\tau]\mapsto T_{[\tau]}$ that is in general not injective. Let $K_{x_0}$ be the kernel of $\rho_M$ and $H_{x_0}=\pi_1(M,x_0)/K_{x_0}$. Since the action of $\pi_1(M,x_0)$ on $\tilde\Sigma$ by $\rho_M$ is isometric, it is already determined by its infinitesimal (orthogonal) action on $T_{x_0}\tilde\Sigma=\nu_{x_0}M$, which is 
$$
\begin{array}{ccc}
\pi_1(M,x_0)\times \nu_{x_0}M&\to& \nu_{x_0}M\\
([\alpha],v)&\mapsto&({\overset{1}{\underset{0}\parallel}}\alpha)v.
\end{array} 
$$
This implies that $H_{x_0}$ is isomorphic to the normal holonomy group of $M$. Thus we can write $\bar M=\tilde M/K_{x_0}$ for the normal holonomy principal bundle $\bar M$, and $H_{x_0}$ is the group of deck transformations of $\tilde M\to\bar M$. \par
Let $\{x_i\}_{i\in I}=M\cap\Sigma$. 
We define $[\sigma_0]=[c_{x_0}]$ and $[\tau_0]=[c_{x_0}]$. For each $i\in I,i\neq 0$, we choose a horizontal curve $[\sigma_i]$ and a vertical curve $[\tau_i]$ from $x_0$ to $x_i$. We can write any element $[\sigma]T_{[\tau]}\in\tilde\Gamma$ as 
$$
\alpha(i,g,h):=g[\sigma_i]T_{h[\tau_i]},
$$ 
where $g\in\pi_1(\Sigma,x_0)$ and $h\in\pi_1(M,x_0)$. 
\begin{lemma}\label{technical}
$\alpha(i,g,h)=\alpha(j,g',h')\iff i=j,  g=g' \mbox{and}\ h^{-1}h'\in K_{x_0}$.
\end{lemma}
\prf $(\Leftarrow)$ follows from $T_{hk[\tau_i]}=T_{[\tau_i]}\circ T_k\circ T_h=T_{[\tau_i]}\circ T_h=T_{h[\tau_i]}$ for $k\in K_{x_0}$, since $T_k=\id_{\tilde\Sigma}$. For $(\Rightarrow)$ apply $[c_{x_0}]$ to both sides. We see that $g[\sigma_i]=\alpha(i,g,h)[c_{x_0}]=\alpha(j,g',h')[c_{x_0}]=g'[\sigma_j]$. The endpoints $x_i$ and $x_j$ are equal, so $i=j$. Thus $g=g'$. It follows $T_h=T_{h'}$, or $T_{h^{-1}h'}=\id_{\tilde\Sigma}$. Observe that $h^{-1}h'\in\pi_1(M,x_0)$. Now we have $h^{-1}h'\in K_{x_0}$.\eop

We remark $\rho_\Sigma(g)=\alpha(0,g,[c_{x_0}])$ for any $g\in\pi_1(\Sigma,x_0)$ and $\rho_M(h)=\alpha(0,[c_{x_0}],h)$ for any $h\in\pi_1(M,x_0)$. \par

\begin{prop}\label{quotient}
The set of leaves $N/\F$ is equal to $\tilde\Sigma/\tilde\Gamma$. If $\Sigma$ is not an exceptional section, $\Gamma$ is defined and preserves $\F$ (specified in the proof). Moreover $N/\F=\Sigma/\Gamma$. The set of sections is $\tilde M/\tilde\Lambda$.
\end{prop}
\prf We have seen at the beginning of this section, that $\hat\pi$ defines a bijection between the set of leaves of $\F$ and that of $\hat \F$. Therefore $N/\F=\hat N/\hat\F$. As justified before we can identify $\Sigma$ with the leaf $T\Sigma$ of $\hat\F^\perp$ in $\hat N$. Under this identification, the section $i_{x_0}:\Sigma\to N$ with image $\Sigma'$ is equal to the restriction $\hat\pi:T\Sigma\to\Sigma'$ of the footpoint map $\hat\pi$ to $T\Sigma$. We also identify the actions of $\tilde\Gamma$ with respect to $\F$ and $\hat\F$. For the bifoliated manifold $\hat N$ we observe that 
$$
\tilde\Gamma[\sigma]=\pi_\Sigma^{-1}(\hat M_{\sigma(1)}\cap T\Sigma)
$$ 
for any $[\sigma]\in \tilde\Sigma$ and hence $\hat N/\hat\F=\widetilde{T\Sigma}/\tilde\Gamma=\tilde\Sigma/\tilde\Gamma$. Assume that $\Sigma$ is not exceptional, so $\Gamma$ is defined. From the formula above we have
$\Gamma(T_x\Sigma)=\hat M_x\cap T\Sigma$ (hence $\hat N/\hat\F=T\Sigma/\Gamma=\Sigma/\Gamma$) and therefore with the obvious identifications $\Gamma(x)=M_x\cap \Sigma$ for regular $x$. We say that $\Gamma$ preserves the restriction $\F_r$ of $\F$ to the regular stratum. The last formula is not well-defined for singular $x\in\Sigma'$, since $\Sigma$ can have self-intersections in $x$, thus there is no canonical element in the preimage of $x$ in $\Sigma=T\Sigma$ under $i_{x_0}$. It is true, however, that any two points $V_1,V_2\in T\Sigma\subset \hat N$ with $\hat\pi(V_i)=x$ lie in the same leaf, namely $\hat M_x$, as seen before. Thus both $\Gamma$-orbits coincide and are equal to $\hat M_x\cap T\Sigma$. Its image under $\hat\pi$ is $M_x\cap \Sigma'$. In this sense $\Gamma$ preserves $\F$, even in singular points. \eop

The description of a $\tilde\Gamma$-orbit in the proof implies in particular that each element of $\tilde\Gamma$ permutes the set $\{g[\sigma_i]\ |\ i\in I, g\in\pi_1(\Sigma,x_0)\}$. In other words, this defines a representation of $\tilde\Gamma$ as a permutation group. This representation is faithful if $M$ has trivial normal holonomy, because of $K_{x_0}=1$ and \Verweisl\ref{technical}.\par
Now let $\F$ be a {\it proper} singular Riemannian foliation admitting sections. Then each regular leaf $M$ has parallel focal structure and finite normal holonomy. The set $\{x_i\}$ is disrete and closed. We call
$$
\CD_{x_i}=\{q\in \Sigma\ |\ d(x_i,q)<d(x_j,q)\ \mbox{for all}\ j\neq i\}
$$
a {\it Dirichlet region} of the set $\{x_i\}$, where $d$ is the distance function in $\Sigma$. These sets are open and disjoint and we have $\bigcup_{i}\overline\CD_{x_i}=\Sigma$. The set $\CD_{x_i}$ is star-shaped and therefore $1$-connected; thus the universal covering $\pi_\Sigma:\tilde\Sigma\to\Sigma$ is trivial over $\CD_{x_i}$ and we denote the connected component of $\pi_\Sigma^{-1}(\tilde\CD_{x_i})$ containing $g[\sigma_i],g\in\pi_1(\Sigma,x_0)$ by $\tilde\CD_{g[\sigma_i]}$. Then $\{\tilde\CD_{g[\sigma_i]}\ |\ g\in\pi_1(\Sigma,x_0),i\in I\}$ is the set of Dirichlet regions for $\pi_\Sigma^{-1}(M\cap\Sigma)$.
\begin{prop} \label{properaction}
Let $\F$ be a proper singular Riemannian foliation admitting sections.
Then the action of $\tilde\Gamma$ on $\tilde\Sigma$ is properly discontinuous. It acts transitively on the set of Dirichlet regions $\{\tilde\CD_{g[\sigma_i]}\ |\ g\in\pi_1(\Sigma,x_0),i\in I\}$ and simply transitive, if $M$ has trivial normal holonomy. The same holds for $\Gamma, \Sigma$ and $\{\CD_{x_i}\}_{i\in I}$ if $\Sigma$ is not an exceptional section. The set of leaves $\tilde\Sigma/\tilde\Gamma$ is an orbifold.
\end{prop}
\prf The action of $\tilde\Gamma$ on $\tilde\Sigma$ is isometric and has discrete orbits, thus it is properly discontinous, i.e., for any compact subset $K$ of $\tilde\Sigma$ the intersection $\phi(K)\cap K$ is non-empty for only a finite number of $\phi\in \tilde\Gamma$. This implies that the set of leaves is an orbifold. The rest follows from \Verweisl\ref{technical}.\eop 

\begin{rem} Singular leaves of $\F$ lift to exceptional leaves. Therefore the nonregular points of the orbifold $\tilde\Sigma/\tilde\Gamma$ correspond exactly to leaves of $\F$ that are either exceptional or singular.
\end{rem}
The following lemma is clear.
\begin{lemma}\label{near_isotropy}
The isotropy group $\tilde\Gamma_{[c_{x_0}]}=\rho_M(\pi_1(M,x_0))\cong H_{x_0}$ is characterized in $\tilde\Gamma$ by mapping $\tilde\CD_{[c_{x_0}]}$ onto itself. Consequently $\tilde\Gamma_{[\sigma]}\subset \tilde\Gamma_{[c_{x_0}]}$ for any $[\sigma]\in\tilde\CD_{[c_{x_0}]}$. An analogous property holds for $\Gamma$ if defined.\eop
\end{lemma}
\begin{rem} Let $G$ be a Riemannian transformation group of $(N,g)$ and let $S$ be a slice through a point $x\in N$ of an orbit $Gx$. It is known that $G_y\subset G_x$ for every $y\in S$. If $Gx$ is an orbit of maximal dimension, this means that the orbit type of $Gx$ is smaller or equal to that of nearby orbits. This corresponds in our theory to $\Gamma_y\subset \Gamma_x$ and that $M_y$ is covering of $M_x$.\par
$H_{x_0}\cong\Gamma_{x_0}$ means that the normal holonomy of a leaf is just the isotropy group of the larger action $\Gamma$. In other words, transversal holonomy generalizes normal holonomy of leaves. 
\end{rem}
We will now give an application for the action of $\tilde\Gamma$. Reinhart showed in \cite{reinhart} that the nearby leaves of a leaf $M$ in a Riemannian foliation are coverings of $M$. The next proposition describes the maximal neighborhood for which this is true. Compare with the proof in \cite{reinhart}.
\begin{prop}\label{nearleafcover}
Let $M$ be a regular leaf of a proper singular Riemannian foliation $\F$ admitting sections and let $x_0\in M$ be arbitrary. Then any regular leaf $M'$ through $\CD_{x_0}$ covers $M$ and the degree is equal to the holonomy orbit $\tilde\Gamma_{[c_{x_0}]}[\gamma]$, where $\gamma$ is a shortest geodesic in $\Sigma$ from $x_0$ to a point in $M'$.
\end{prop}

\prf Let $y_0\in \CD_{x_0}\cap M'$ and let $\gamma_0$ be a shortest geodesic from $x_0$ to $y_0$ which is contained in $\CD_{x_0}$. 
Moreover let $Y:=\{[\gamma_j]\}_{j\in J}:=\tilde\Gamma_{[c_{x_0}]}[\gamma_0]$. We define an action $h\cdot ([\tau],[\gamma_j])=(h[\tau],T_{[h^{-1}]}[\gamma_j])$ of $\pi_1(M,x_0)$ on $\tilde M\times Y$. Note that this group acts from the left by \Verweisl\ref{recrules}, and the action is free and properly discontinuous. Let $\tilde M\times_{\pi_1(M,x_0)} Y:=(\tilde M\times Y)/\pi_1(M,x_0)$
We want to show that
$$
\begin{array}{ccc}
\tilde M\times_{\pi_1(M,x_0)}Y&\to& M'\\ \\
([\tau],[\gamma_j])&\mapsto&(T_\tau\gamma_j)(1)
\end{array}  
$$
is a diffeomorphism. The map is clearly surjective. We show that it is well-defined. Let $h\in\pi_1(M,x_0)$ and $([\tau],[\gamma_j])\in\tilde M\times Y$. Then
$$
(T_{T_{h^{-1}}\gamma_j}h\tau)(1)=(T_{T_{h^{-1}}\gamma_j}h\cdot 
T_{T_hT_{h^{-1}}\gamma_j}\tau)(1)=(T_{\gamma_j}\tau)(1),
$$
so the map is well-defined. We prove injectivity. Let $(T_{\gamma_j}\tau)(1)=(T_{\gamma_k}\tau')(1)$ for vertical curves $\tau,\tau'$ starting at $x_0$. Then $\tau(1)=\tau'(1)$, so there is exactly one $h\in \pi_1(M,x_0)$ such that $h[\tau]=[\tau']$. We claim $\gamma_k=T_{h^{-1}}\gamma_j$. We have $(T_{\tau}\gamma_j)(1)=(T_{\gamma_j}\tau)(1)=(T_{\gamma_k}\tau')(1)=(T_{\tau'}\gamma_k)(1)$. Thus
$(T_\tau\gamma_j)(1)=(T_{h\tau}\gamma_k)(1)=(T_\tau(T_h\gamma_k))(1)$. Applying $T_{\tau^{-1}}$ shows
$\gamma_j(1)=(T_h\gamma_k)(1)$, i.e., $\pi_\Sigma([\gamma_j])=\pi_\Sigma(T_h[\gamma_k])$. Since $[\gamma_j]$ and $T_h[\gamma_k]$ lie in $\tilde\CD_{[c_{x_0}]}$ this implies $[\gamma_j]=T_h[\gamma_k]$ and we proved our claim. Now the above map is a diffeomorphism. Thus $M'=\tilde M\times_{\pi_1(M,x_0)}Y$ covers $M$ with typical fiber $Y$. \eop

With the introduced methods we have extended \Verweisp \ref{necessary_foliation} in \cite{toeben}:
\begin{prop}\label{parallel-cutlocus}
Let $M$ be a closed and embedded submanifold with parallel focal structure and finite normal holonomy. If $v\in\nu M$ is a multiplicity 
$k$ focal normal of vertical type so are its normal parallel translations. In other words the vertical focal data is also invariant under normal parallel translation. If $v$ is a cut normal, so are its normal parallel 
translations. In particular, the cut distance function is constant along the parallel normal 
fields. \eop
\end{prop}
The last statement is already known from \Verweisp \ref{necessary_foliation}. With this proposition we can distinguish a submanifold with parallel focal structure from other submanifolds by its cut locus. The next proposition proved in \cite{toeben} shows the relation between $\{\CD_{x_i}\}$ and the cut locus of a regular leaf, which has parallel focal structure as we know. It also shows that a leaf $M'$ through $\CD_{x_0}$ is regular if $M$ has a globally flat normal bundle.
\begin{prop}\label{fundamentaldomain-cutlocus}
$\bigcup_{i\in I}\partial{\CD_{x_i}}\subset \CC_{(M,N)}$ and $\CC_{(M,N)}\cap \CD_{x_i}$ are points of the cut locus of $\exp^\Sigma_{x_i}$. If $M$ has a globally flat normal bundle, then the focal points of horizontal type are contained in $\bigcup_{i\in I}\partial{\CD_{x_i}}$.\eop
\end{prop}
Now we express \Verweisp\ref{nearleafcover} as a corollary in terms of the cut locus.
\begin{corollary}\label{cover-cutlocus}
Let $M$ be a closed and embedded submanifold with parallel focal structure and finite normal holonomy. Then the parallel submanifolds that are not contained in the cut locus of $M$ are coverings of $M$.\eop
\end{corollary}

The following result is a another corollary of \Verweisp\ref{nearleafcover}.
\begin{corollary}\label{trivialhol}
The regular leaves with a globally 
flat normal bundle are diffeomorphic to each other. They cover any other regular leaf, the exceptional leaves. These exceptional leaves, if they exist, are contained in the cut locus of any regular leaf with a globally flat normal bundle. The union of regular leaves with a globally flat normal bundle is open and dense in $N$.\eop
\end{corollary}

The exceptional leaves lie in the cut locus of the leaves with trivial normal holonomy. Can we give a sufficient condition on $N$ that guarantees that all regular leaves have trivial normal holonomy? We will show that there are no exceptional leaves, if the ambient space is a simply connected symmetric space. We need some preparations. For a point $p\in N$ we define $\CP=\CP(N,\varphi\times p)$ as the set of pairs $(x,\gamma)$, where $x\in M$ and $\gamma:[0,1]\to N$ is a $H^1$-curve in $N$ with $\gamma(0)=\varphi(x)$ and $\gamma(1)=p$. We write $\CP(N,M\times p)$ for the path space if $\varphi:M\to N$ is the inclusion map. It is known that $\CP$ is a Hilbert manifold. The smooth function 
$$
E_p\, :\, \CP\to\RR;\ (x,\gamma)\mapsto\int_0^1\|\dot\gamma(t)\|^2\,dt
$$
is called the {\it energy functional} (associated to $p$). The map $E_p$ is a Morse function, i.e., it has only non-degenerate critical points, if and only if $p$ is not a focal point of $\varphi$. We assume that $p$ is not a focal point, i.e., $p$ is regular for the normal exponential map of $M$. The energy functional is bounded below by zero and it is known that it satisfies the Palais-Smale condition. For $s\in\RR$ we write $\CP^s=E_p^{-1}\{[0,s]\}$ and $\CP^{s-}=E_p^{-1}\{[0,s)\}$. Let $s$ be a regular value of $E_p$. The Morse inequalities state $b_k(\CP^s)\leq \mu_k(E_p|\CP^s)$, where $b_k(\CP^s)$ is the $k$-th Betti number of $\CP^s$ with respect to $\ZZ_2$ and $\mu_k(E_p|\CP^s)$ is the number of critical points of index $k$ of $E_p$ below $s$. 
\begin{thm}\label{symmetric}
Let $\varphi:M\to N$ be a proper immersion with parallel focal structure with finite normal 
holonomy into a simply connected symmetric space $N$.
Then $\varphi$ factorizes finitely over an embedding 
that has a globally flat normal bundle. There are no
exceptional parallel submanifolds and the cut locus of $M$ only consists of focal points.
\end{thm}
The first part of the next result was proven in \Verweisl 1A.3 of \cite{podestathorbergsson} for polar actions. The proof refers to \Verweisl 2.10 in \cite{ewert}, which has a gap. We want to present the complete proof and relate the occurrence of exceptional leaves to the cut locus.
It is well-known that a closed hypersurface $M$ of a simply connected manifold $N$ is orientable and thus has a globally flat normal bundle. If the codimension is greater than one we have to argue differently.\medskip

\prf Since the image of $\varphi$ is a leaf and because $\varphi$ is proper, it factorizes finitely over an embedding. So we can assume that $\varphi$ is this embedding. Let $p$ be a regular point in $N$ with respect to the normal exponential map of $M$. We claim that $E_p$ has only one (local) minimum. We will prove this later. Assume now that there is an exceptional parallel manifold $M'$ of $M$. We choose $\eps>0$ smaller than the injectivity radius of $M'$, which is positive since the cut distance function is constant, $p\in M'$ and $v\in\nu_pM'$ with non-trivial holonomy degree and $\|v\|< \eps$. Let $w\neq v$ in $\nu_p M'$ be a normal parallel translation of $v$. Let $M:=M_v'$. The geodesics $\gamma_v|[0,1]$ and $\gamma_w|[0,1]$, if parameterized in reverse direction, are normal to $M$, nonfocal and of index $0$ by the choice of $\eps$. Therefore there are at least two minima of $E_p:\CP(N,M\times p)\to\RR$, contradiction. Now assume that there is point $p$ in the cut locus of $M$ that is not a focal point. By \Verweisp \ref{focal-cut-point} it is the endpoint of two minimal normal geodesics of index $0$, contradiction.\par
We will now prove the claim. First we observe that $\CP$ is connected because of the homotopy sequence of the fibration $\CP\to M;c\mapsto c(0)$ and $\pi_1(N)=1$. Let $\gamma\in\CP$ be an arbitrary critical point of index 1 of $E_p$, 
i.e., $\gamma$ is a normal geodesic of index 1, with $E_p(\gamma)=\kappa$. Let $\eps>0$ so small that there is no critical level in $[\kappa-\eps,\kappa+\eps]$ other than $\kappa$. Let $e^1$ be the corresponding $1$-cell in $\CP^{\kappa+\eps}$ through $\gamma$ attached to $\CP^{\kappa-\eps}$. Let $v\in\nu_x M, x=\gamma(0)$ with $\gamma_v=\gamma$, and let $t_0v$, $ 0<t_0<1$ be the focal normal with multiplicity 1. First we assume that $v$ is a focal normal of horizontal type. Let $F$ be the focal leaf of $v$ in $M$ (see \Verweisr \ref{focalleaf}) and not in $\bar M$. Since $F$ is $1$-dimensional and compact we have $F\cong S^1$. We construct a variation $\lambda:F\to\CP$ of $\gamma$ by
$$
\lambda(y)(t):=
\left\{\begin{array}{ll}
\eta(t v_y)&\mbox{if}\ t\in [0,t_0]\\
\gamma(t)&\mbox{if}\ t\in [t_0,1]
\end{array}\right.
$$
(compare with \cite{thorbergssondupin} and \cite{ewert}).
This smooth map is injective (if we took the focal leaf in $\bar M$ this map would only be a covering) and Ewert deforms it under the negative gradient flow of $E_p$ to a map $\lambda':F\to\CP$ that has a unique non-degenerate maximum in $x$. 
We denote the generator of $H_1(F)=H_1(S^1)$ by $[S^1]$, where we consider homology over $\ZZ_2$. Then $z:=\lambda_*([S^1])=\lambda'_*([S^1])\in H_1(\CP^{\kappa+\eps})$ is a so-called {\it Bott-Samelson cycle} with the property $j_{1}(z)=[e^1]$, where $j_1:H_1(\CP^{\kappa+\eps})\to H_1(\CP^{\kappa+\eps},\CP^{\kappa-\eps})$ is the map of the homology sequence of the pair $(\CP^{\kappa+\eps},\CP^{\kappa-\eps})$. Now we assume that $v$ is a focal normal of vertical type, i.e., $\gamma(t_0)$ is conjugate to $x$ along $\gamma$ in $\Sigma_x$ with multiplicity $1$. Since $\Sigma_x$ is a symmetric space as a totally geodesic submanifold of $N$, an $S^1$-action fixing $x$ and $\gamma(t_0)$ applied to $\gamma|[0,t_0]$ gives an $S^1$-familiy of geodesics from $\gamma(0)$ to $\gamma(t_0)$. We extend this variation as above to a map $\lambda:S^1\to\CP$ and Ewert proves that also $z:=\lambda_*([S^1])\in H_1(\CP^{\kappa+\eps})$ is a Bott-Samelson cycle with $j_1(z)=[e^1]$. Now let us consider a part of the homology sequence of the pair $(\CP^{\kappa+\eps},\CP^{\kappa-\eps})$:
$$
H_{1}(\CP^{\kappa+\eps})\stackrel{j_1}{\longrightarrow}  H_{1}(\CP^{\kappa+\eps},\CP^{\kappa-\eps})\stackrel{\partial_{1}}{\longrightarrow} H_{0}(\CP^{\kappa-\eps})\stackrel{i_0}{\longrightarrow} H_{0}(\CP^{\kappa+\eps})
$$
By standard Morse theory the set of $[e^1]$ for the descending cells $e^1$ of all critical points $\gamma$ of index $1$ on level $\kappa$ is a basis of $H_{1}(\CP^{\kappa+\eps},\CP^{\kappa-\eps})$. As seen above, each $[e^1]$ is in the image of $j_1$. Therefore $\partial_1=0$. Thus
$i_0:H_0(\CP^{\kappa-\eps})\to H_0(\CP^{\kappa+\eps})$ is injective for every critical level $\kappa$ and consequently also $i_0:H_0(\CP^r)\to H_0(\CP)$ for every $r$. Since $\CP$ is connected, $H_0(\CP)=\ZZ_2=H_0(\CP^r)$, so $\CP^r$ is connected for every $r$. If $E_p$ had at least two minima, the higher one say on level $\kappa$, $\CP^{\kappa+\eps}$ would be disconnected for small $\eps$, contradiction.\eop

Note that we can drop symmetry of $N$ if we assume the sections are either symmetric or do not have conjugate points.

\end{document}